\PassOptionsToPackage{table,dvipsnames}{xcolor}
\documentclass[12pt]{article}

\usepackage[english]{babel}
\usepackage[utf8]{inputenc}
\usepackage{geometry}
\usepackage{multicol,float}
\usepackage{multirow}
\usepackage{array}
\usepackage[affil-it]{authblk}

\usepackage{amsthm}
\usepackage{amsmath}
\usepackage{amsfonts}
\usepackage{amssymb}
\usepackage{mathtools}  
\usepackage{graphicx}
\usepackage{todonotes}
\usepackage{complexity}
\usepackage{adjustbox}
\usepackage{changepage}

\usepackage{algorithm,algpseudocode}
\usepackage{longtable}

\usepackage{cancel}
\usepackage{comment}
\usepackage[inline]{enumitem}
\usepackage{ulem}
\usepackage{url}
\usepackage{subcaption}
\usepackage{empheq}
\usepackage{appendix}
\usepackage{natbib}
\setcitestyle{authoryear,open={(},close={)}} 

\usepackage{hyperref}
\hypersetup{
	colorlinks   = true, 
	urlcolor     = cyan, 
	linkcolor    = blue, 
	citecolor   = blue 
}

\makeatletter
\def\blfootnote{\xdef\@thefnmark{}\@footnotetext}
\makeatother

\makeatletter
\def\ps@pprintTitle{%
  \let\@oddhead\@empty
  \let\@evenhead\@empty
  \let\@oddfoot\@empty
  \let\@evenfoot\@oddfoot
}

\definecolor{light-gray}{gray}{0.95}

\newcommand{\norm}[1]{\| #1 \|}

\title{\textbf{Aircraft conflict resolution: A benchmark generator}\blfootnote{\textit{The current affiliations of the authors are:}\\
\textbf{Mercedes Pelegr\'in}: FICO Xpress Optimization (Global Analytics Delivery); \textit{email address:} \texttt{mariamercedes.pelegrin@gmail.com} \\
\textbf{Martina Cerulli}*: Department of Computer Science, University of Salerno, 84084 Fisciano, Italy; \textit{email address:} \texttt{mcerulli@unisa.it} }}
\date{}

\author[1]{Mercedes Pelegr\'in}
\author[2]{Martina Cerulli\thanks{Corresponding author}}

	\affil[1]{\small Laboratoire d'Informatique de l'X (LIX), \'Ecole Polytechnique,  91128 Palaiseau, France}
	
	\affil[2]{\small Information Systems, Decision Sciences and Statistics Department, ESSEC Business School of Paris, 95000 Cergy-Pontoise, France}

\begin{document}
\maketitle
\vspace*{-1cm}
\hspace*{-6mm}\fcolorbox{red}{white}{\parbox{\textwidth}{This paper has been accepted for publication in the \textbf{INFORMS Journal on Computing}. The \textbf{final published version} is available at \url{https://doi.org/10.1287/ijoc.2022.1265}. The software that supports the findings of this study is available at the IJOC GitHub software repository (\url{https://github.com/INFORMSJoC/2021.0283}).}}
\vspace*{4mm}

\begin{abstract}
Aircraft conflict resolution is one of the major tasks of computer-aided Air Traffic Management and represents a challenging optimization problem.
Many models and methods have been proposed to assist trajectory regulation to avoid conflicts. 
However, the question of testing the different mathematical optimization approaches against each other is still open. Standard benchmarks include unrealistic scenarios in which all the flights move towards a common point or completely random generated instances.
There is a lack of a common set of test instances that allows comparison of the available methods under a variety of heterogeneous and representative scenarios.
We present a flight deconfliction benchmark generator that allows the user to choose between 
\begin{enumerate*}
\item[(i)] different predefined scenario inspired on existing benchmarks in the literature;
\item[(ii)] pseudo-random traffic meeting certain congestion measurements;
\item[(iii)] and randomly generated traffic.
\end{enumerate*} 
The proposed setting can account for different levels of  difficulty in the deconfliction of the aircraft and allows to explore and compare the real limitations of optimization approaches for aircraft conflict resolution.

The software that supports the findings of this study is available within the paper and its Supplemental Information (\url{https://pubsonline.informs.org/doi/suppl/10.1287/ijoc.2022.1265}) as well as from the IJOC GitHub software repository (\url{https://github.com/INFORMSJoC/2021.0283}) at (\url{http://dx.doi.org/10.5281/zenodo.7377734}).
\end{abstract}

\section{Introduction}
Decision support tools for Air Traffic Management (ATM) has been on the focus of latest developments in airspace systems,
including American NextGen and European SESAR projects.
Commercial aviation, together with emerging new actors in the context of urban air mobility, need for computational tools that enable some level of automatization in airspace operations to have scalable systems.
This need has driven the efforts of the Computer Science research community, including the Operations Research community, towards answering the challenge of automatization of ATM. 
The main goal has been to give support to human air traffic controllers. In this domain, an air sector is monitored by the controller, who needs to provide separation provisions about 20-30 minutes before conflicts can occur. Mathematical models typically receive a set of rectilinear trajectories crossing the observed air sector as an input; the task is to provide new trajectories that are conflict-free. This problem is known as Conflict Resolution or Tactical Deconfliction (TD), as opposed to strategic deconfliction, which happens hours before flight, and collision avoidance, which addresses imminent conflicts.

The TD problem can be stated as follows: given a set of $n$ aircraft with initial positions $\hat p_i$, $i=1,\ldots,n$, and nominal vectors of velocity $\hat V_i$, $i=1,\ldots,n$, find new vectors $V_i$ such that the resulting trajectories are conflict-free, that is, such that
\begin{equation}\label{eq:mindist}
\norm{(\hat p_i+tV_i)-(\hat p_j+tV_j)}\geq D \qquad \forall\; i,j=1,\ldots,n, \forall \;t\geq 0,
\end{equation}
where $D$ is the minimum safety distance. The new vectors correspond to adjustments made on the nominal ones, namely heading angles and/or speed changes. The aircraft trajectories can be considered in the Euclidean 2D space (when a fixed altitude is assumed) or in the 3D space. 
We refer the interested reader to \citet{survey-separation} for a detailed analysis of equation \eqref{eq:mindist} and equivalent representations of aircraft separation. In particular, it is sufficient that \eqref{eq:mindist} is satisfied at the time instant where the minimum distance between $i$ and $j$ is attained. We use this last condition in our code implementation.

Computational testing of mathematical models for TD, including performance comparison, needs from reference sets of meaningful benchmarks.
In this work, we present a TD instances generator (available in GitHub \citep{TDinstancesGenerator}) that provides a common source for evaluation for both 2D and 3D cases.
To the best of our knowledge, no TD instances generator existed in the literature so far, even though standard benchmarks have been collected in some repositories, see e.g. \cite{repRey,repMarti,repWang} (presented in \cite{cerulli,dias21,wang} respectively).

Standard TD benchmarks can be classified into two main groups, namely those based on predefined scenarios and those generated randomly. Scenario-based instances represent synthetic, usually unrealistic, configurations that are used to stress conflict resolution methods. 
Such instances, including the so-called circle, rhomboidal and grid instances, have been mostly used for testing approaches in 2D. Here, we standardize the generation of these benchmarking instances, providing a reference tool for testing. Moreover, we extend these predefined scenarios to the 3D case.

Randomly generated instances (with no control on the density of conflicts) have been used in previous literature, e.g., in \cite{escudero16,caf-rey17,frazzoli}. In most of the cases, such instances are generated by placing a set of aircraft in a $2$-dimensional airspace represented by a square, with random speeds and heading angles. However, the lack of a common reference makes the different studies not comparable. We propose a pseudo-random traffic generator, where traffic congestion is determined by the user through certain input parameters. This allows to generate TD instances with a predefined complexity, something that, to the best of our knowledge, was not proposed before. Finally, the generator can also produce completely random instances, as the ones already used in literature, in a given rectangular region. 

A benchmark instance is defined by the initial positions $\hat p_i$ and nominal vectors $\hat V_i$. As output of the generator, several additional information regarding the generated instance are given, including number of conflicts, and, for each pair of aircraft in conflict, duration of the conflict and pairwise distance at the time of minimal separation (see the toy examples provided in the repository \cite{TDinstancesGenerator} with illustrative purposes).
With this information we give an idea of the air traffic complexity \citep{sylvie,isufaj} of the generated instances. The most widely used metric for aircraft complexity is indeed aircraft density, which measures the number of aircraft flying within an air sector. This metric, however, is shown insufficient to capture Air Traffic Controllers (ATCs) workload, as studied in \cite{delahaye}. In this sense, the generator we propose displays a varied set of complexity indicators for the generated instances. In the case of our pseudo-random instances, the tunable complexity indicators are related to traffic density. That is, the user can decide on the  conflicts density  of the generated instances, but they cannot control e.g. the duration of the conflicts. 

The rest of the paper is organized as follows.
The scenario-based instances, namely, circle, rhomboidal, grid, and their extensions to 3D (resp. sphere, polyhedral, and cubic) are described in Section \ref{sec:scenario}. Section \ref{sec:random} is devoted to the pseudo-random and random instances. Section \ref{sec:compu} present a computational study to demonstrate the accuracy of the proposed method in generating pseudo-random benchmarks. Finally, Section~\ref{sec:conclusion} closes the paper with some concluding remarks and potential future extensions. A list of the parameters used in the paper can be found in the Appendix. The commands that can be used to generate the instances illustrated in the different sections of the paper may be found in the repository \cite{TDinstancesGenerator}.

\section{Scenario-based instances}\label{sec:scenario}
In this section we describe the six predefined scenarios available in our generator, which account for the 2D and 3D implementations of three different layouts that exist in the literature. In each scenario, both the initial positions and nominal vectors of velocity of the aircraft are generated according to a predefined configuration. We provide a random variant for each of the six scenarios, which consists in adding a random deviation to the predefined nominal vectors of velocity. The range of such random deviation can be selected by the user. Each scenario is characterized by certain parameters, which can be tuned by the user. For further technical details on the options and parameters of the generator, see the repository \cite{TDinstancesGenerator}. 

\subsection*{Circle problem}
\begin{figure}[h!]
\begin{subfigure}[t]{0.45\textwidth}
 	\begin{center}
 	\includegraphics[scale=0.3]{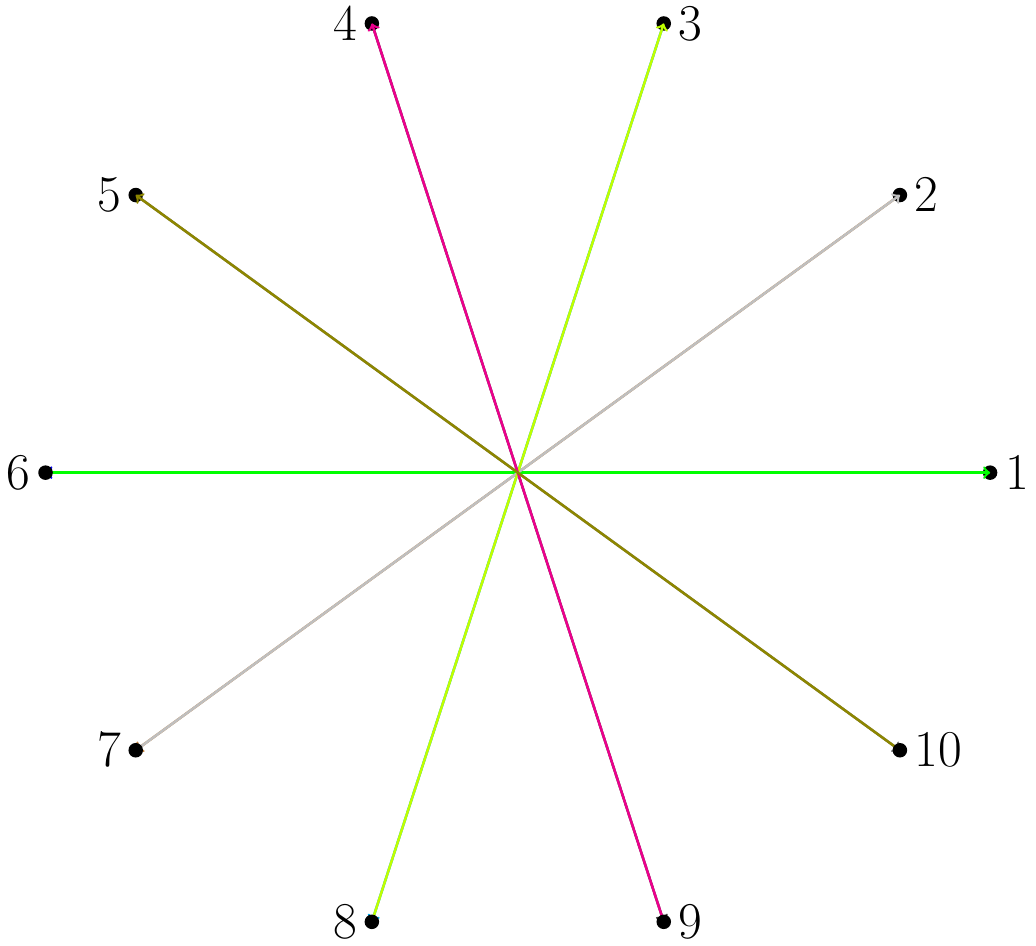}
	\end{center}
	\caption{Standard circle}
	\label{fig:CP-a}
\end{subfigure}
\hspace{0.5cm}
\begin{subfigure}[t]{0.45\textwidth}
 	\begin{center}
 	\includegraphics[scale=0.4]{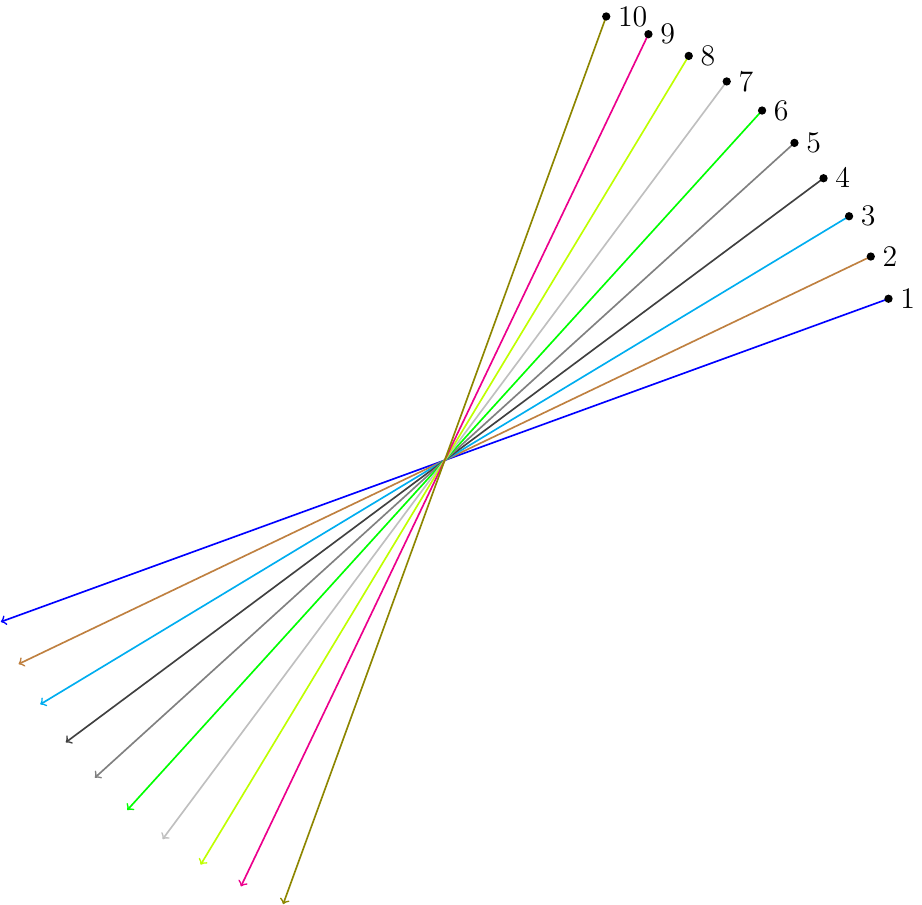}
	\end{center}
	\caption{Sector circle}
	\label{fig:CP-b}
\end{subfigure}
\caption{Circle problem with $n=10$}
\label{fig:CP}
\end{figure}

\begin{figure}[h!]
\begin{subfigure}[t]{0.45\textwidth}
 	\begin{center}
 	\includegraphics[scale=0.65]{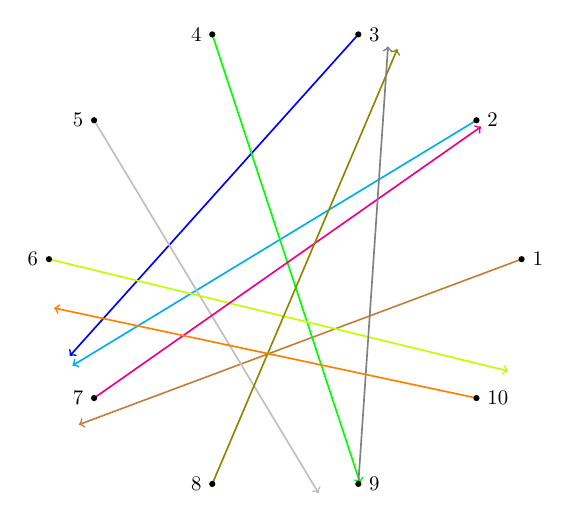}
	\end{center}
	\caption{Standard random circle}
	\label{fig:RCP-a}
\end{subfigure}
\hspace{0.5cm}
\begin{subfigure}[t]{0.45\textwidth}
 	\begin{center}
 	\includegraphics[scale=0.65]{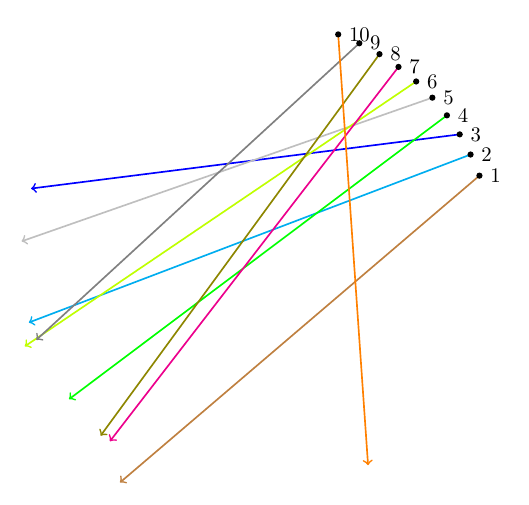}
	\end{center}
	\caption{Sector random circle}
	\label{fig:RCP-b}
\end{subfigure}
\caption{Random circle problem with $n=10$}
\label{fig:RCP}
\end{figure}

The circle problem refers to a scenario in which $\hat p_i$  are uniformly distributed on a circle, and $\hat V_i$  point towards its center. These instances (presented for the first time in \cite{frazzoli}, under the name ``symmetric encounter patterns"), although very unrealistic, have served to stress existing resolution methods in order to explore their limits. Indeed, assuming uniform speed for all the aircraft, the number of conflicts of these kind of instances is maximum, $n(n-1)/2$,  all of them taking place at the same time instant. Benchmarking circle instances may be found in \cite{repRey}.
 
This scenario is characterized by the number of aircraft, $n$, and the radius of the circle. Given these two parameters, the aircraft initial positions are uniformly distributed on the circumference and their trajectories point towards the center, see Figure~\ref{fig:CP-a}.
Alternatively, one can select a sector of the circumference where the initial positions are uniformly distributed, see Figure~\ref{fig:RCP-b}. The beginning of the sector and its width are input parameters that can be tuned by the user. The default option is that of the sector equal to the whole circle (standard circle problem). This ``sector option" was inspired by some works \citep{cafieri,rey16}, where aircraft are distributed on a quarter of the circumference. 

The so-called random circle problem, firstly presented in \cite{vanaret}, and then used in \cite{reyrapine,rey16,caf-rey17} among others, were thought to provide more realistic scenarios by randomly deviating the trajectories of the aircraft from the center of the circle.
Figure~\ref{fig:RCP} illustrates two random circle instances. Figures~\ref{fig:RCP-a} and \ref{fig:RCP-b} respectively correspond to the standard and sector options, and we can observe that the aircraft trajectories do not necessarily pass through the center of the circle in this case. 
Benchmarking random circle instances are contained in the repositories \cite{repRey}, and \cite{repWang}. As far as it concerns \cite{repWang}, aircraft at different altitude level are considered, and also the initial position of the aircraft is randomly chosen within a square having a center lying on the circumference.

\subsection*{Sphere problem}
\begin{figure}[h!]
\begin{subfigure}[t]{0.45\textwidth}
 	\begin{center}
 	\includegraphics[scale=0.4]{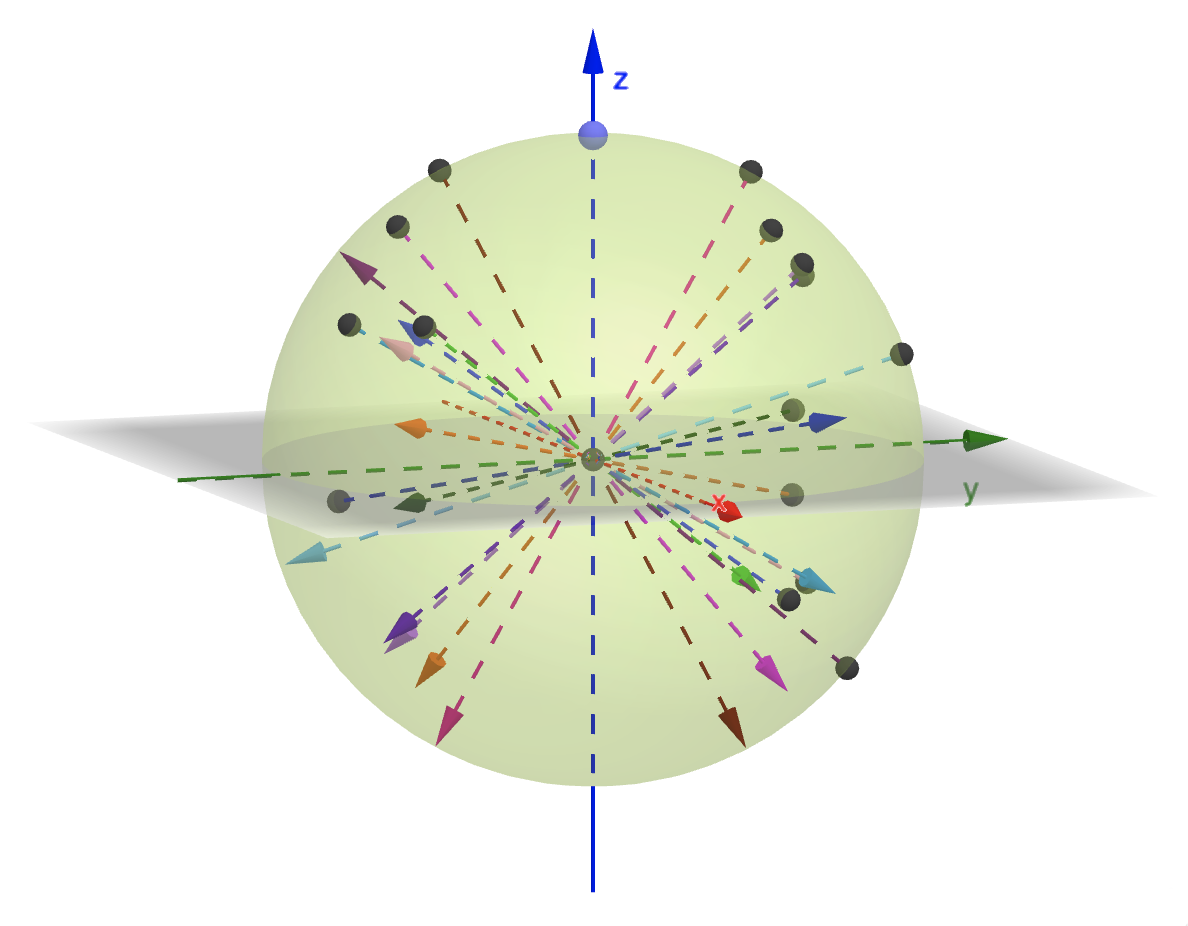}
	\end{center}
	\caption{Standard sphere}
	\label{fig:sphere-a}
\end{subfigure}
\hspace{0.5cm}
\begin{subfigure}[t]{0.45\textwidth}
 	\begin{center}
 	\includegraphics[scale=0.4]{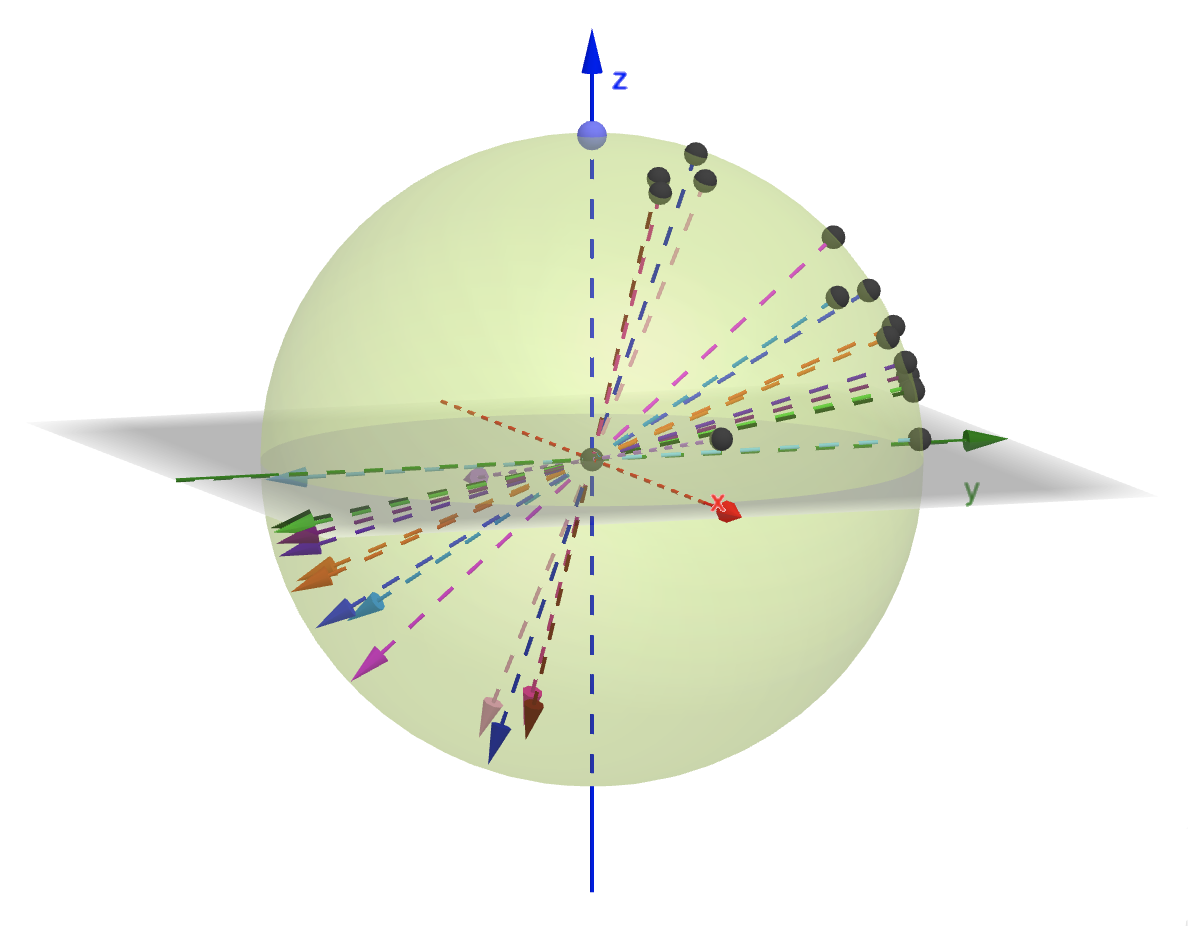}
	\end{center}
	\caption{Sector sphere}
	\label{fig:sphere-b}
\end{subfigure}
\caption{Sphere problem with $n=15$ }
\label{fig:sphere}
\end{figure}

We extend the circle problem to the {3D} case, which we name sphere instance, see Figure~\ref{fig:sphere}\footnote{All 3D figures in this paper have been generated using GeoGebra \citep{geogebra}}. Instances of these type may be found in the repository \cite{repMarti}. In the sphere instances, $n$ aircraft are uniformly distributed on the surface of a sphere of a given radius~$r$. 

\begin{figure}[h!]
\begin{subfigure}[t]{0.45\textwidth}
 	\begin{center}
 	\includegraphics[scale=0.4]{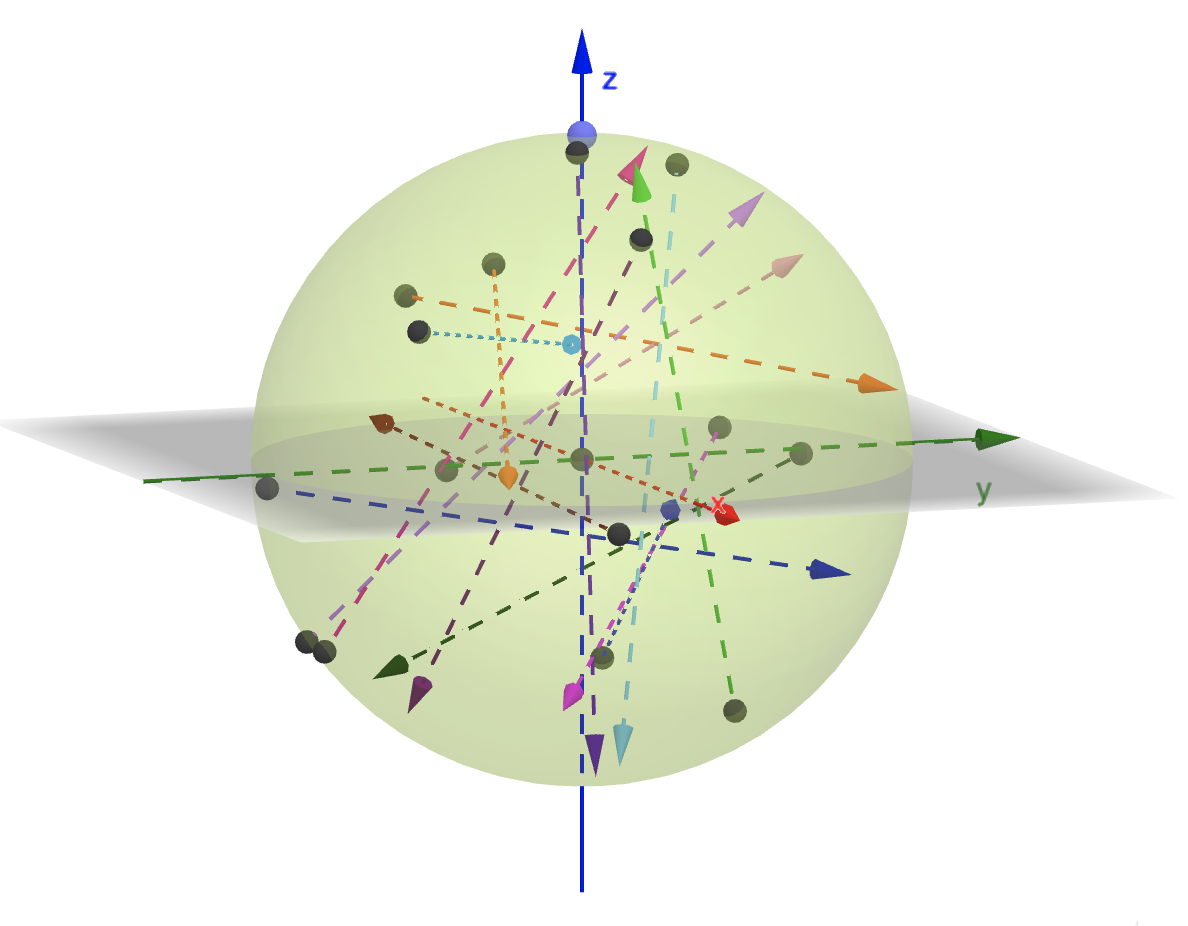}
	\end{center}
	\caption{Standard random sphere}
	\label{fig:random-sphere-a}
\end{subfigure}
\hspace{0.5cm}
\begin{subfigure}[t]{0.45\textwidth}
 	\begin{center}
 	\includegraphics[scale=0.4]{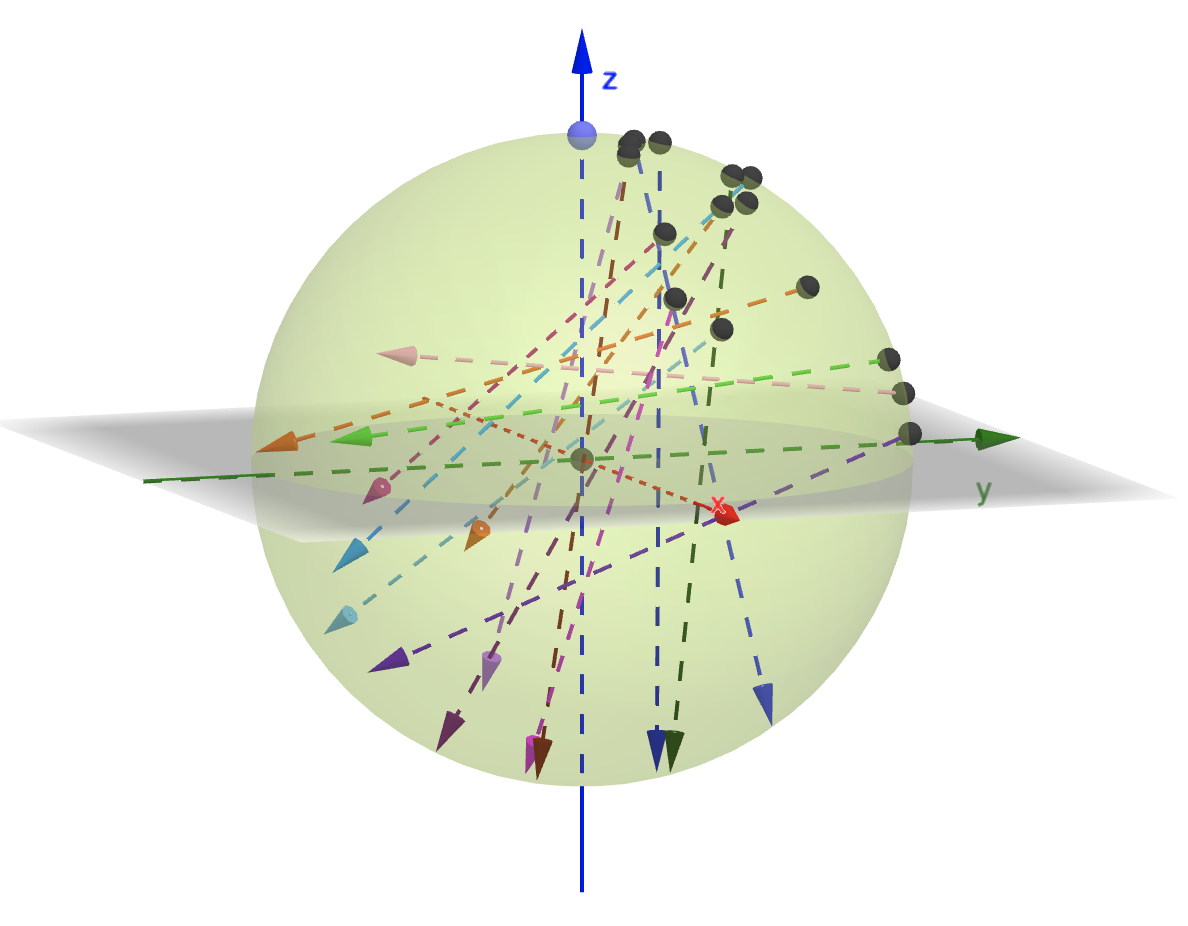}
	\end{center}
	\caption{Sector random sphere}
	\label{fig:random-sphere-b}
\end{subfigure}
\caption{Random sphere problem with $n=15$}
\label{fig:random-sphere}
\end{figure}
Let $\hat \theta_i$ and $\hat \phi_i$ denote the heading angles that the trajectory of aircraft $i$ forms with the $x$-axis and $z$-axis respectively. We generate angles $\hat \theta_i,\hat \phi_i$, $i=1\ldots n$, in such a way that, when applying the following transformation into Cartesian coordinates
$$\hat{x_i}=r\sin{\hat \phi_i}\cos{\hat \theta_i}, \quad \hat{y_i}=r\sin{\hat \phi_i}\sin{\hat \theta_i}, \quad \hat{z_i}=r\cos{\hat \phi_i}, $$ we obtain a uniform distribution of points $\hat{p_i}=(\hat{x_i},\hat{y_i},\hat{z_i})$ on the sphere.
The so-called random sphere problem is characterized, similarly to the random circle problems, by aircraft trajectories slightly deviating w.r.t. the center of the sphere, see Figure~\ref{fig:random-sphere}.

The ``sector" option is also available (see Figures~\ref{fig:sphere-b} and \ref{fig:random-sphere-b} for an example).

\subsection*{Rhomboidal problem} 
In this scenario, aircraft fly on predefined rectilinear trails. The trails are organized in two groups. The trails within each of the groups are parallel, while every pair of trails of different groups cross at a given angle, $\alpha$. Without loss of generality, one group is made of horizontal trails (with slope equal to 0) and the other one of slopping trails (with slope equal to $\alpha$). 

\begin{figure}[H]
	\begin{subfigure}[t]{0.3\textwidth}
		\begin{center}
			\includegraphics[scale=0.7]{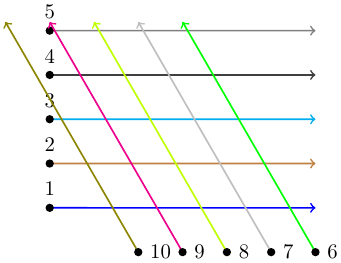}
		\end{center}
		\caption{$\alpha=120$ degrees}
		\label{fig:rhomboidal-a}
	\end{subfigure}
	\begin{subfigure}[t]{0.34\textwidth}
		\begin{center}
			\includegraphics[scale=0.7]{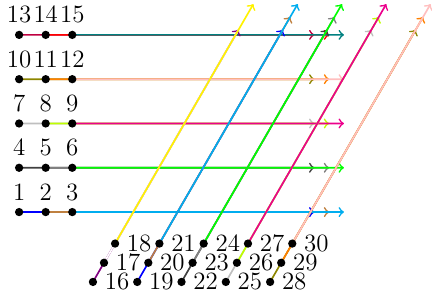}
		\end{center}
		\caption{3 flights per trail}
		\label{fig:rhomboidal-b}
	\end{subfigure}
	\begin{subfigure}[t]{0.4\textwidth}
		\begin{center}
			\includegraphics[scale=0.7]{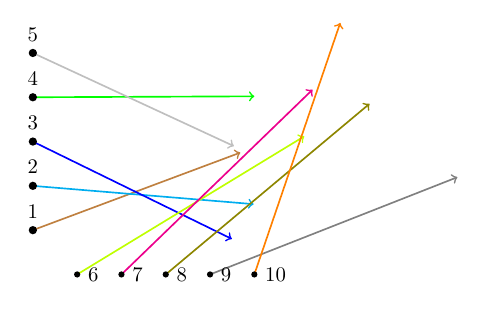}
		\end{center}
		\caption{Random rhomboidal poblem}
		\label{fig:rhomboidal-c}
	\end{subfigure}
	\caption{Rhomboidal problem with $10$ trails}
	\label{fig:rhomboidal}
\end{figure}

For instance, Figure~\eqref{fig:rhomboidal-a} shows a rhomboidal scenario with 10 trails (5 horizontal and 5 slopping) and $\alpha=120^o$. This example features one aircraft per trail, i.e., $n=10$. However, one can generate more than one aircraft flying on the same trail, as depicted in Figure~\ref{fig:rhomboidal-b}. The separation between consecutive aircraft on the same trail can be tuned by the user, as well as the separation between parallel trails. The random option is also available for the rhomboidal scenario, and is illustrated in Figure~\ref{fig:rhomboidal-c}. In this case, the trails are deviated by an angle randomly chosen in an interval, which can be set by the user.

The number of aircraft in the rhomboidal scenario is fully characterized by the specific configuration of the trails. That is, if we call $m_x$ and $m_y$ the number of horizontal and slopping trails, and $n_x$ and $n_y$ the number of aircraft at each horizontal and slopping trail, respectively, we have
$n= m_x\cdot n_x+m_y\cdot n_y.$
Then, $n$ should not be input by the user for the rhomboidal scenario. 

The rhomboidal (or flow) scenario was introduced in \cite{frazzoli} with the name ``Crossing aircraft stream", and instances of this type have been used as benchmarks in \cite{lehouillier,cafieri}. 
Also, they can be interpreted as a generalization of the flow instances available at \cite{repRey}, where only two trails are considered. These kind of instances provide a more realistic scenario than that of the circle ones, but still ensuring a certain level of congestion. Namely, assuming uniform speed of the aircraft and uniform spacing of the trails, the number of conflicts is equal to the number of intersections between the trails times the minimum between $n_x$ and $n_y$.

\subsection*{Polyhedral problem}

We define the polyhedral instances as a generalization of the rhomboidal ones for the 3D case (basic polyhedral instances may be found in the repository \cite{repMarti}).
The idea is to have several 2D horizontal planes at different altitudes (different $z$-components). Each of these planes contains trails following a rhomboidal layout. In addition, vertical planes containing slopping trails of aircraft are also included to complete the configuration of the scenario.

\begin{figure}[H]
\begin{subfigure}[t]{0.5\textwidth}\hspace{-1.5cm}
\begin{center}
\includegraphics[scale=0.4]{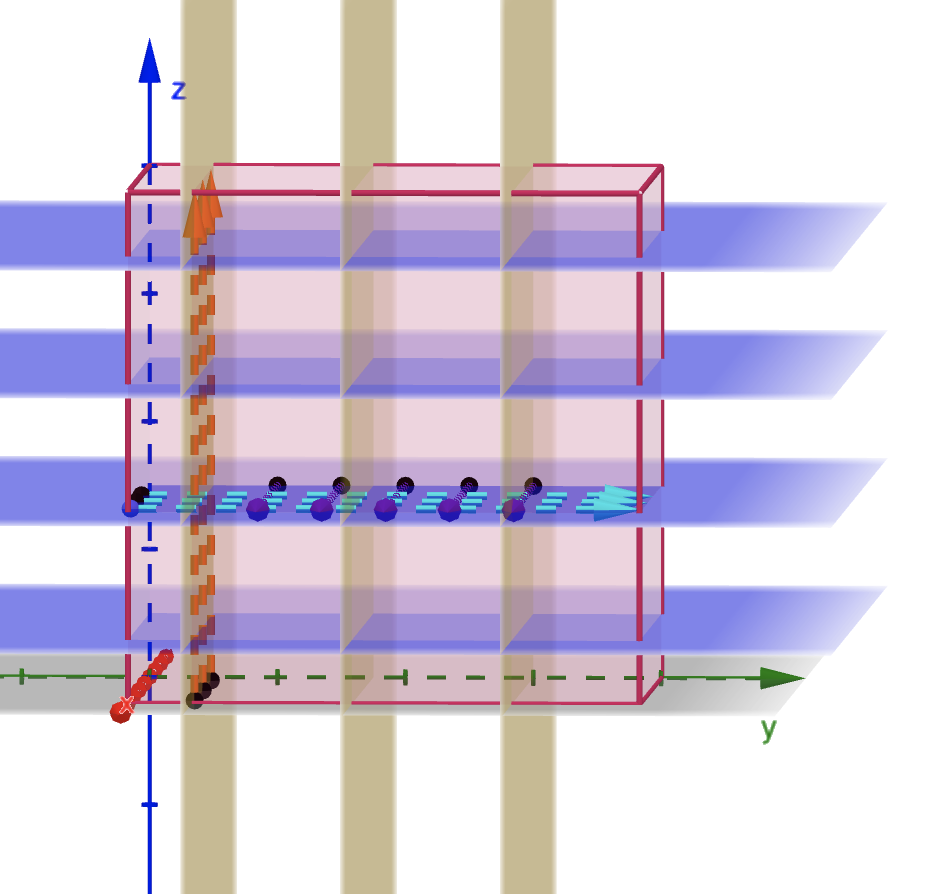}
\end{center}
\caption{Planes where trails are generated}
\label{fig:3Dplanes}
\end{subfigure}\hspace*{-2.2cm}
\begin{subfigure}[t]{0.5\textwidth}
 	\begin{center}
 	\includegraphics[scale=0.4]{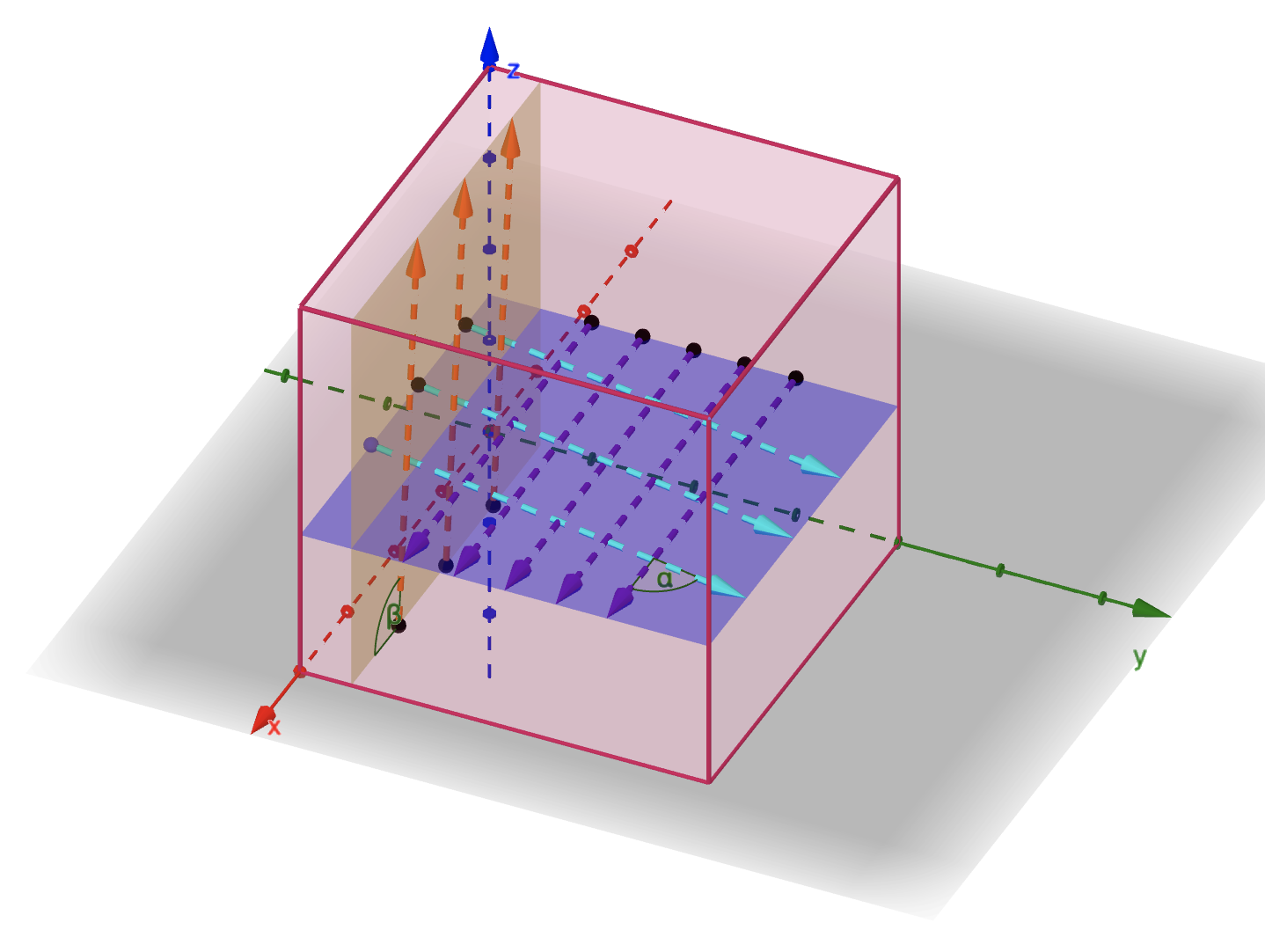}
	\end{center}\vspace*{-1cm}
	\caption{Polyhedral instance}
	\label{fig:poly}
\end{subfigure}
\begin{subfigure}[t]{1\textwidth}
 	\begin{center}
    \includegraphics[scale=0.2]{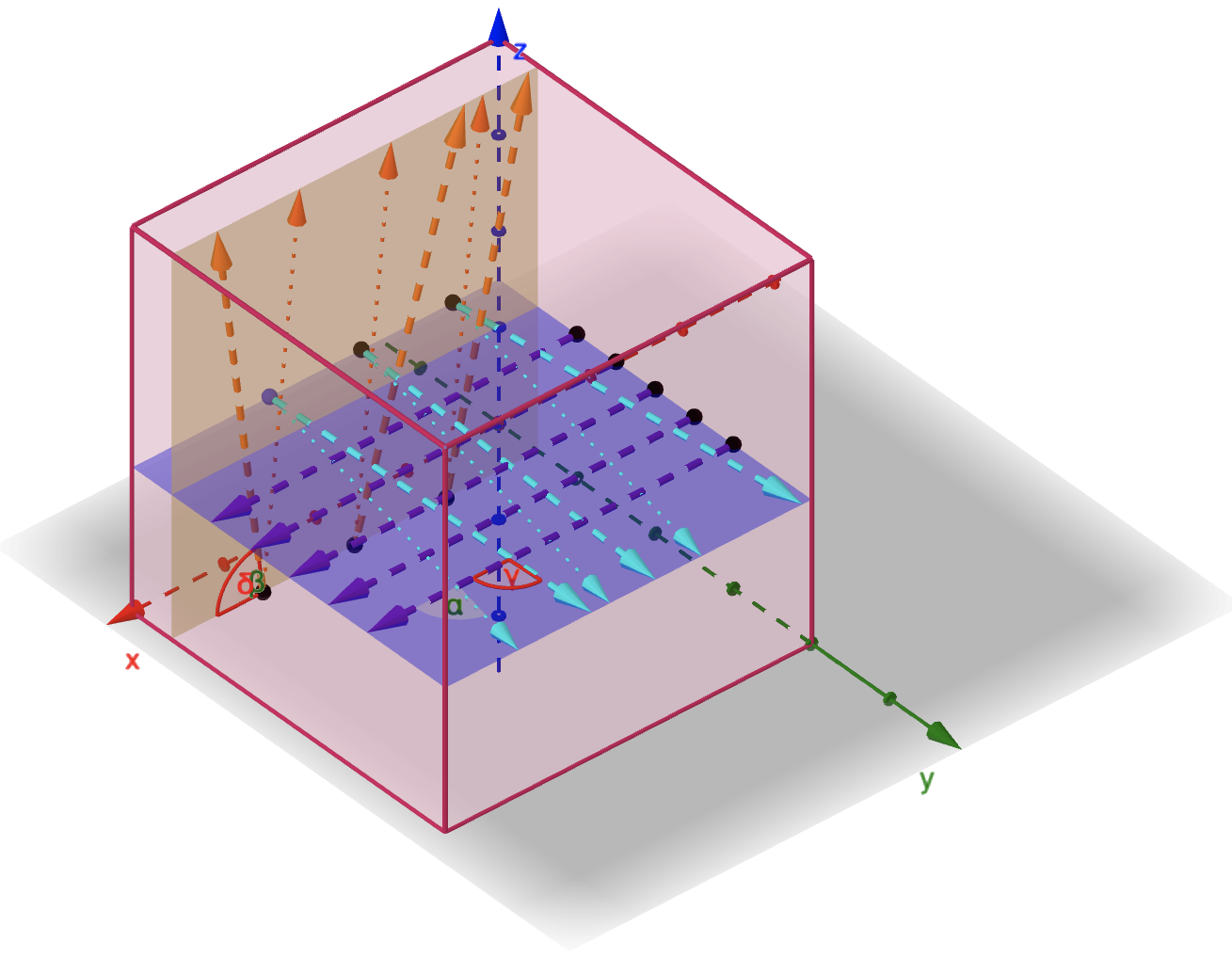}
    \end{center}\vspace{-0.8cm}
    \caption{Random polyhedral instance}
    \label{fig:my_label}\label{fig:random-poly}
\end{subfigure}
\caption{Illustration of 3D polyhedral scenario}
\label{fig:3D}
\end{figure}

Figure~\ref{fig:3D} illustrates polyhedral instances. In Figure~\ref{fig:3Dplanes}, we see the family of horizontal parallel planes in blue, and that of vertical parallel planes in yellow, which intersect the former. Both the number of horizontal (resp. vertical) planes, and the distance among them can be tuned by the user.
Figure~\ref{fig:poly} illustrates the layout of the trails inside the planes. 
Each blue plane can be seen as a 2D rhomboidal scenario: it contains a set of trails that are parallel to $x$-axis (they would correspond to the ``horizontal trails" in the 2D case) and another one with slope $\alpha^h$, which can be set different at each horizontal plane $h$. On the other hand, each yellow plane contains a set of parallel trails that cross the horizontal planes forming an angle that we call $\beta^v$, which, again, may be set different at each vertical plane $v$. As in the 2D case, there can be one or more aircraft flying on each trail. The separation between consecutive aircraft on the same trail, as well as the number of aircraft on each trail can be tuned by the user.

In this scenario, the user can specify a different number of horizontal/slopping trails for each horizontal plane, as well as a different number of slopping trails for each vertical plane. On the other hand, the distance between parallel trails can be decided only for all horizontal/vertical planes, and not for each of them individually.

We also consider the random polyhedral problem, illustrated in Figure~\ref{fig:random-poly}, where both $\alpha^h$ and $\beta^v$ can be deviated by an angle randomly chosen in an interval, the boundaries of which can be given by the user.

\subsection*{Grid problem} 
\begin{figure}[H]
\begin{subfigure}[t]{0.32\textwidth}
 	\begin{center}
 	\includegraphics[scale=0.7]{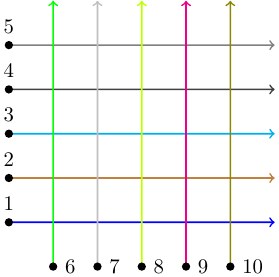}
	\end{center}
	\caption{1 flight per trail}
	\label{fig:grid-a}
\end{subfigure}
\begin{subfigure}[t]{0.32\textwidth}
 	\begin{center}
 	\includegraphics[scale=0.7]{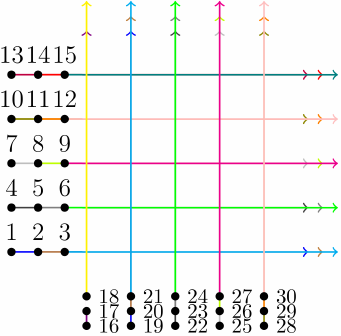}
	\end{center}
	\caption{3 flights per trail}
	\label{fig:grid-b}
\end{subfigure}
\begin{subfigure}[t]{0.32\textwidth}
 	\begin{center}
 	\includegraphics[scale=0.7]{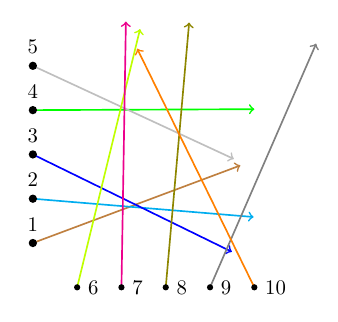}
	\end{center}
	\caption{Random grid}
	\label{fig:grid-c}
\end{subfigure}
\caption{Grid problem with $10$ trails}
\label{fig:grid}
\end{figure}
The grid instances are particular cases of the rhomboidal scenario with $\alpha=\pi/2$. 
Figure~\ref{fig:grid} illustrates grid configurations with one aircraft per trail (Figure~\ref{fig:grid-a}), three  aircraft per trail (Figure~\ref{fig:grid-b}), and a random grid instance (Figure~\ref{fig:grid-c}). Note that random grid and random rhomboidal instances (Figure~\ref{fig:rhomboidal-c}) look alike. They differ in the original angle $\alpha$ to which a random deviation is applied.

The grid scenario was firstly presented in \cite{niedringhaus}. We can find a grid instance among the benchmarks used in \cite{cafieri}. The repository \cite{repRey} also includes grid instances, although these are in fact particular cases of the ones presented in this work (they feature only 4 trails, 2 horizontal and 2 vertical).

\subsection*{Cubic problem}

\begin{figure}[H]\hspace*{-2cm}
	\begin{subfigure}[t]{0.6\textwidth}
		\begin{center}
			\includegraphics[scale=0.4]{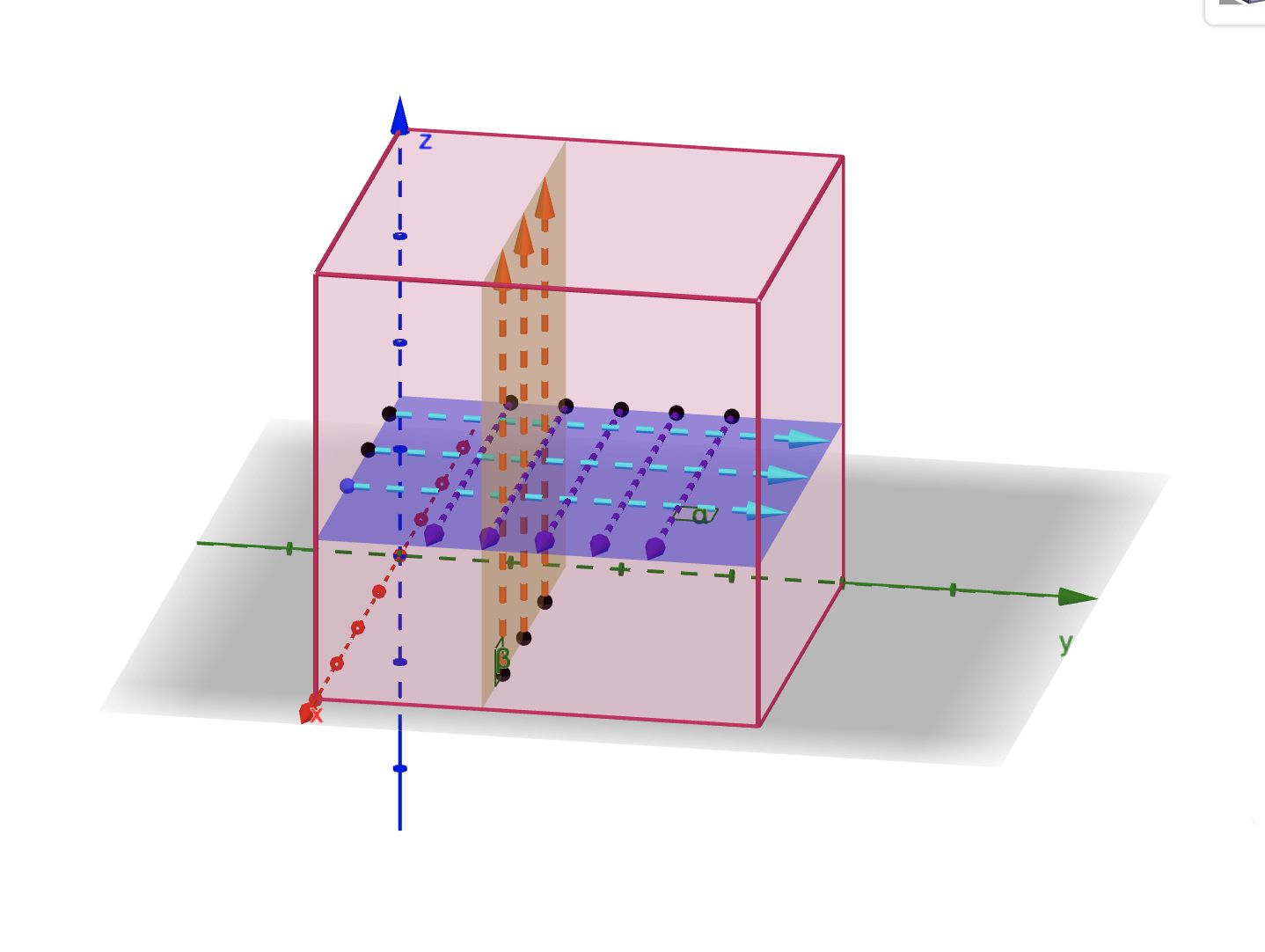}\vspace*{-1.5cm}
		\end{center}
		\caption{Cubic problem}
		\label{fig:3Dgrid-a}
	\end{subfigure}
	\begin{subfigure}[t]{0.5\textwidth}
		\begin{center}
			\includegraphics[scale=0.4]{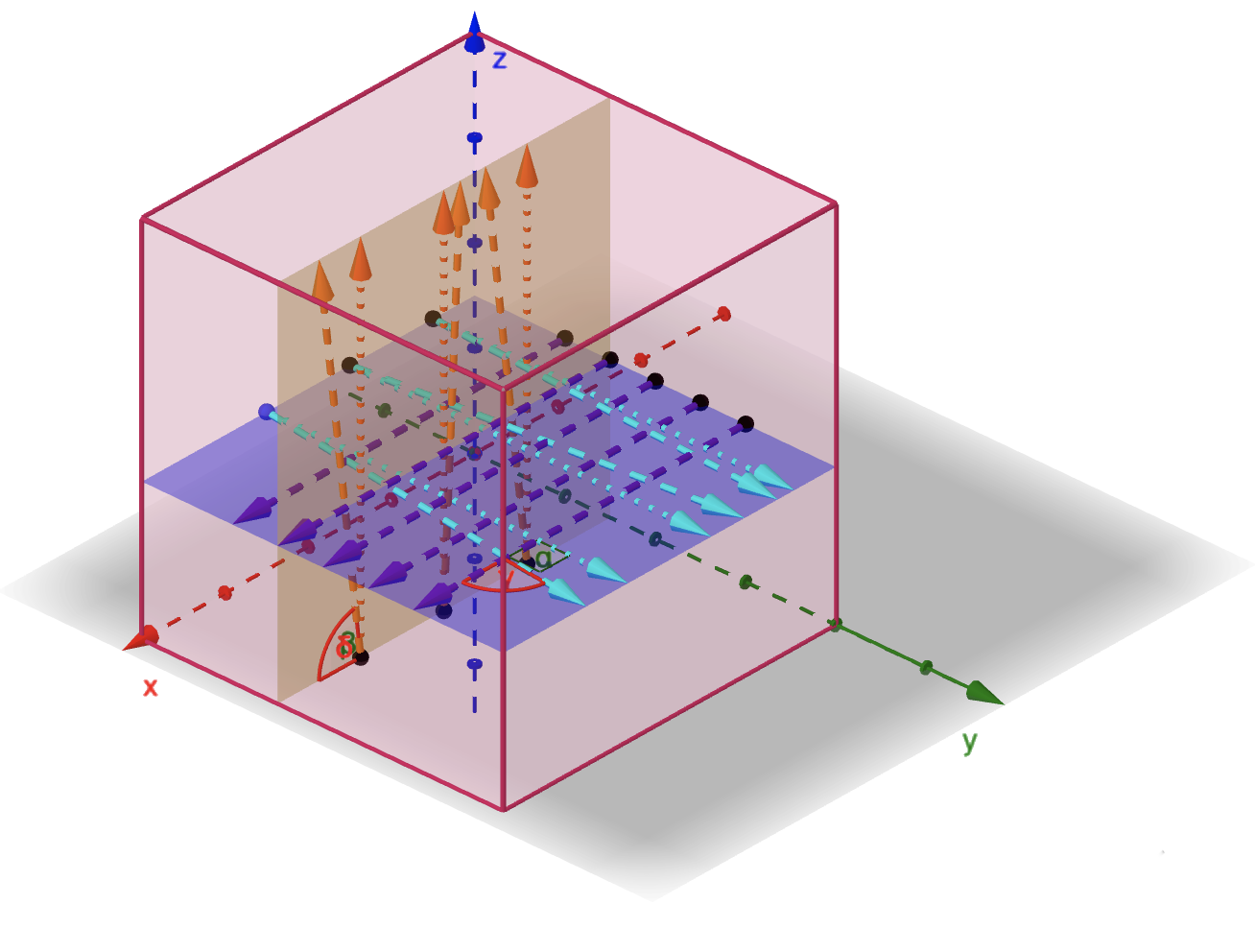}\vspace*{-1cm}
		\end{center}
		\caption{Random cubic problem}
		\label{fig:3Dgrid-b}
	\end{subfigure}
	\caption{Illustration of the 3D cubic scenario}
	\label{fig:3Dgrid}
\end{figure}

To generalize the grid scenario to the 3D case, we define the cubic instances. These are particular cases of polyhedral instances with $\alpha^h=\beta^v=\pi/2$ for all $h$ and $v$ (see Figure~\ref{fig:3Dgrid}, where Figure~\ref{fig:3Dgrid-b} shows the random case). We note that the cubic problem can be particularly representative for Urban Air Mobility (UAM) settings, where vertical take-off and landing aircraft (VTOLs) come into play. 

\section{Random and pseudo-random instances}\label{sec:random}

In contrast to scenario-based random instances presented in the previous section, we describe here random instances where the aircraft trajectories are not obtained deviating by a certain angle from a specific layout. Instead, we consider a rectangle/parallelepiped (air sector) in which the initial positions $\hat p_i$ and nominal vectors $\hat V_i$ are generated randomly or pseudo-randomly.

On the one hand, we provide the option of randomly generating both 2D and 3D instances in order to offer such common reference. On the other hand, we address pseudo-random instances for the first time (to the best of our knowledge).
This consists in generating aircraft configurations that meet a particular level of congestion, which is determined by the user through different parameters.
Namely, we define the following parameters related to  traffic congestion:
\begin{itemize}
\item $n_c$: total number of pairs of aircraft that are in conflict;
\item $p_c$: probability that one aircraft is in conflict with at least another aircraft;
\item $max_c$: maximum number of aircraft that a fixed aircraft can be in conflict with.
\end{itemize}
The expected number of conflicts is given by:
\begin{equation}\label{eq:nconf}
     E(n_c)= n \cdot p_c \cdot \left(\frac{1+max_c}{2}\right)\cdot \frac{1}{2}.
\end{equation}
Therefore, the user would only need to introduce two of these three parameters (the remaining one would be calculated according to \eqref{eq:nconf}).
The maximum number of conflicts that can be generated is given by $(n\cdot max_c/2)$, thus we discourage the use of values of  $n_c$ and $max_c$ such that $max_c< 2\cdot n_c/n$. On the other hand, the average number of conflict is at most $(n\cdot (max_c+1)/4)$ (if $p_c=1$). Then, if both $n_c$ and $max_c$ are input by the user, we recommend to use values such that $max_c\geq  4\cdot n_c/n-1$. The user may decide to input just some of these three parameters or none of them. If no parameter is specified, then $p_c=0.5$, $max_c=n-1$ by default, and $n_c$ is fixed to the closer integer to $E(n_c)$ (see the repository \cite{TDinstancesGenerator} for more details).

In the following, we describe the algorithm used to generate a given set of trajectories meeting the desired level of congestion. Our implementation is made of two main parts. First, the initial positions are generated, $\hat p_i$, for all $i=1,\ldots,n$. 
We generate these positions near the borders of a given rectangle (or near the faces of a given parallelepiped in the 3D case), and the heading angles pointing towards its interior.

\begin{figure}[H]
 	\begin{center}
 	\includegraphics[scale=0.4]{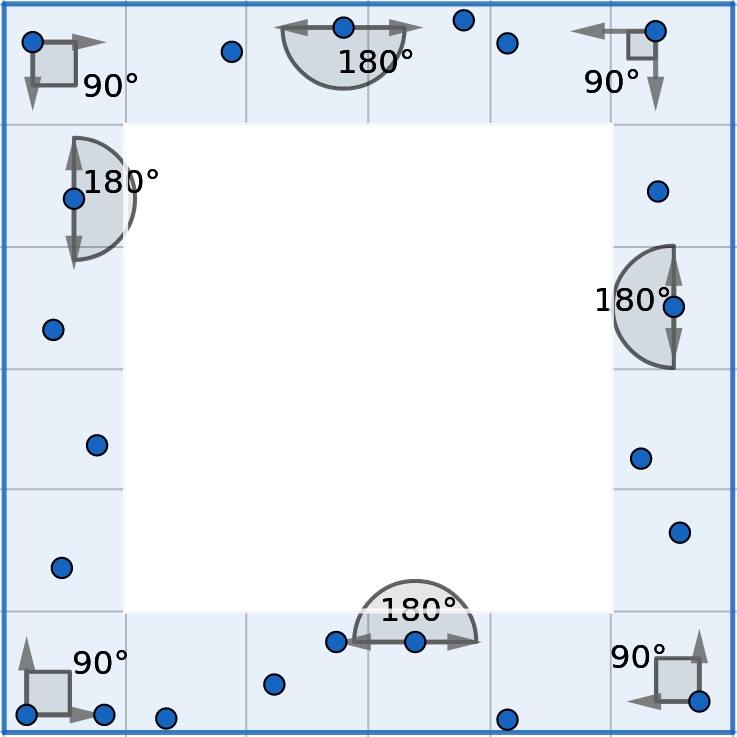}
	\end{center}
	\caption{Initial positions of aircraft in a 2D pseudo-random instance}
	\label{fig:inipos2d}
\end{figure}

The idea is that the observed air sector (the rectangle or the parallelepiped) was empty at a previous time instant, and several aircraft are about to cross the sector now. Figure~\ref{fig:inipos2d} illustrates the initial positions in a pseudo-random 2D instance, which are randomly generated inside squared areas of a given width near the borders.

\begin{figure}[H]
 	\begin{center}\vspace{-0.5cm}
 	\includegraphics[scale=0.4]{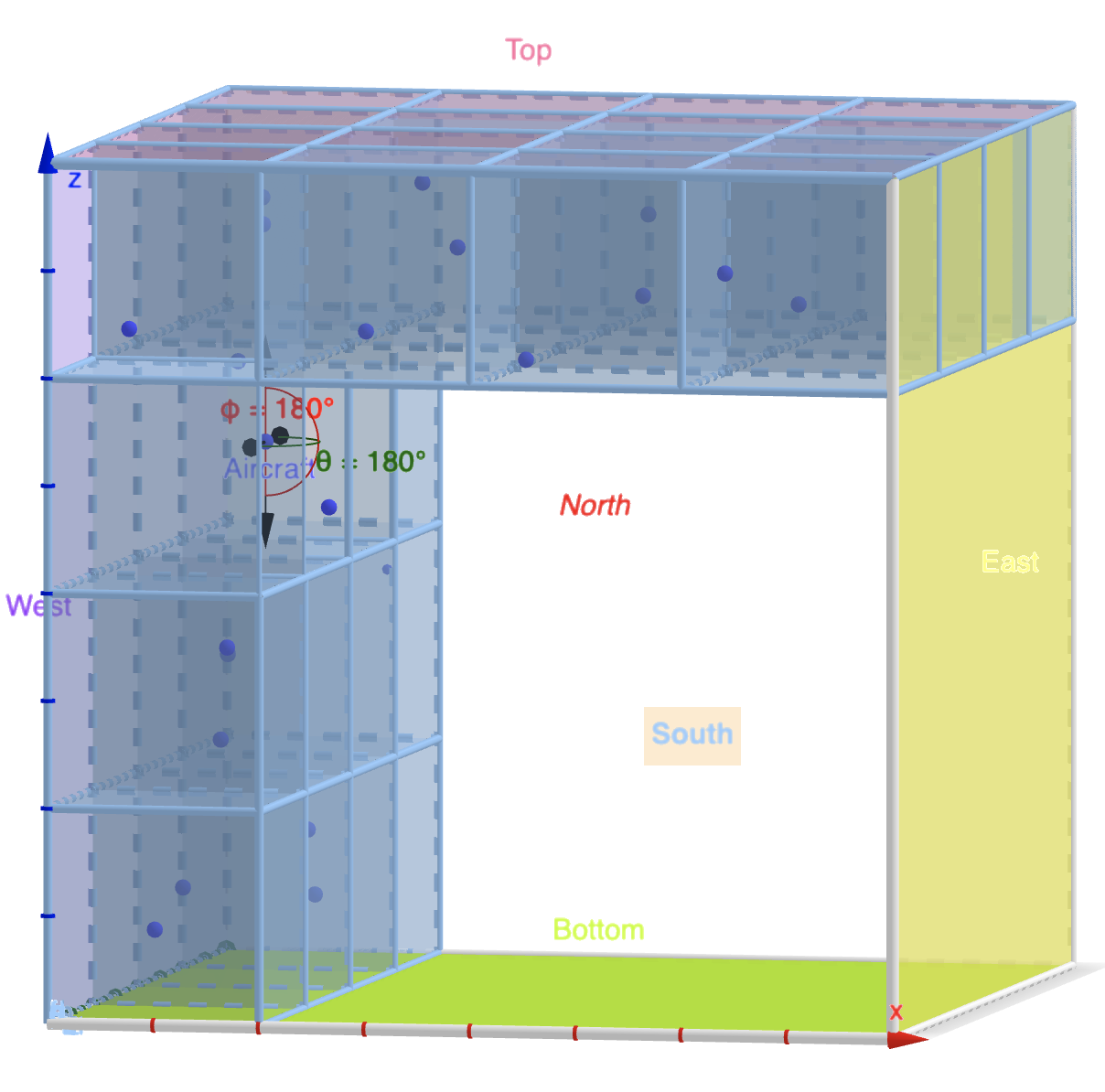}
	\end{center}
	\caption{Initial positions of aircraft in a 3D pseudo-random \textit{west-top} scenario}
	\label{fig:inipos3d}
\end{figure}

Figure~\ref{fig:inipos3d} shows the initial positions of aircraft in a pseudo-random 3D instance, randomly generated inside cubic areas of a given volume near west and top faces of the air sector. 

The generator allows the initial positions to be near the 4 borders of a rectangle (for 2D instances) or the 6 faces of a parallelepiped (3D case), but also permits other configurations as that depicted in Figure~\ref{fig:inipos3d}. In particular, we identify the following orientation of the borders/faces to define our configurations:
\begin{itemize}
    \item \textit{north} and \textit{south} borders/faces as the ones normal to the vector of the $y$-axis (with the south one adjacent to the hyperplane $y\equiv 0$),
    \item \textit{west} and \textit{east} borders/faces as the ones normal to the vector of the $x$-axis (with the west one adjacent to the hyperplane $x\equiv 0$),
    \item (in the 3D configuration) \textit{top} and \textit{bottom} faces as the ones normal to the vector of the $z$-axis (with the bottom one adjacent to the hyperplane $z\equiv 0$).
\end{itemize}
The possible additional configurations are characterized by aircraft only on one of the following pairs of borders of the rectangle for the 2D instances: north/south, west/east, or west/north. For the 3D instances, the following pairs of faces of the parallelepiped characterize the allowed additional configurations: north/south, west/east, west/north, north/top, west/top (see Fig.~\ref{fig:inipos3d}), or top/bottom. As depicted in Figures \ref{fig:inipos2d} and \ref{fig:inipos3d}, each aircraft is assigned valid heading angles ranges depending on its position.

After all $\hat p_i$ are generated, the vectors $\hat V_i$ have to be added with the aim of attaining the target congestion. This process is described in Algorithm~\ref{alg1}. At each step of the main while loop, a different aircraft $i$, $i\in\{1,\ldots,n\}$, is randomly selected. For this aircraft, a target number of conflicts ($targetConf$ in Algorithm~\ref{alg1}) between 1 and $max_c$ is randomly generated with probability $p_c$ (see Step 1.4). After generating $\hat V_i$, the resulting conflicts between $i$ and the aircraft processed in previous iterations are counted. If they are equal to the target number of conflicts $targetConf$, we continue with the next iteration of the main while loop. That is, a new randomly chosen aircraft will be processed. Otherwise, $maxTrials$ defines the maximum number of times a new vector of velocity $\hat V_i$ is generated in order to obtain the target number of conflicts, $targetConf$. If this maximum is achieved, other distinct target numbers of conflicts are tried up to a certain limit, defined by $maxConf$. Lastly, if none of the values of target conflicts in $\{0,\ldots,maxConf\}$ was successful, higher values would be tried starting from $maxConf+1$ and increasingly adding one unit.

The pseudocode in Algorithm \ref{alg1} describes our implementation. The input/output parameters are defined above and in Table~\ref{tab:parameters}. The variables updated throughout the algorithm are: 
\begin{itemize}
    \item $isExplored[i]$: true if $\hat V_i$ already has been assigned a value, $i=1,\ldots,n$
    \item $nExplored$: number of aircraft for which  $\hat V_i$  has been generated
    \item  $totalConf$: number of conflicts generated until the moment.
\end{itemize}
The function $randInt(a,b)$ used in Step 1.4 generates an integer number in the range $[a,b]$. 
The function $countConf(i,isExplored)$ used in Step 1.6 returns the number of conflicts between $i$ and all aircraft $j$ s.t. $isExplored[j]=true$.

While we generate the trajectories, in every iteration of the main loop in Step 1, $p_c$ is updated in Step 1.11 so that it is consistent with the already generated conflicts. In particular, it is set to $4\frac{n_c-totalConf}{(n-nExplored)*(1+max_c)}$, according to Eq.~\eqref{eq:nconf} when dealing with the remaining $(n-nExplored)$ aircraft for which we want to generate the remaining $(n_c-totalConf)$ conflicts.

\begin{figure}[H]
\vspace*{-0.5cm}
\begin{subfigure}[t]{0.5\textwidth}
 	\begin{center}
 	\includegraphics[scale=0.5]{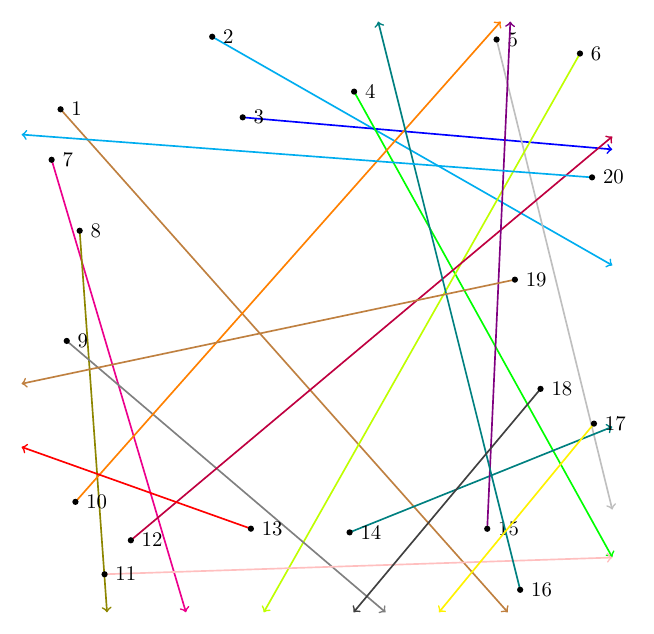}
	\end{center}
	\caption{{\footnotesize2D pseudo-random problem with $n=20$, $n_c=21$}}
	\label{fig:randEx}
\end{subfigure}
\hspace{0.3cm}
\begin{subfigure}[t]{0.6\textwidth}
 	\begin{center}
 	\includegraphics[scale=0.6]{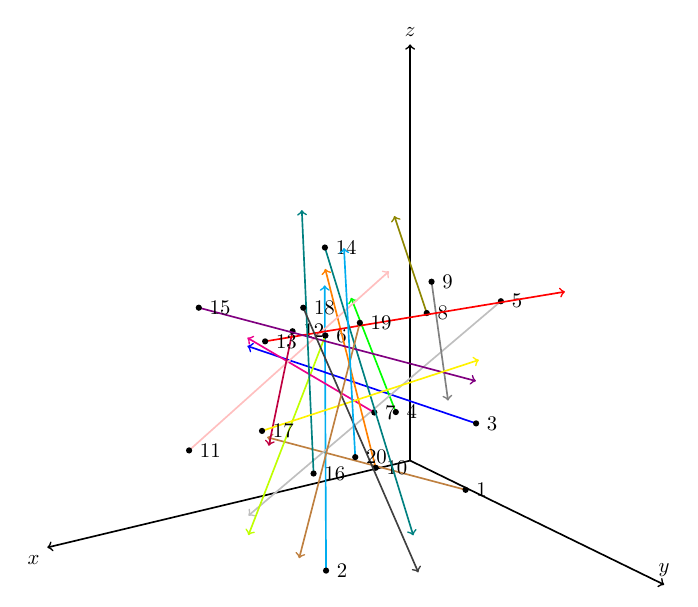}
	\end{center}
	\caption{{\footnotesize3D pseudo-random problem with $n=20$, $n_c=13$}}
	\label{fig:randEx3D}
\end{subfigure}
\caption{Pseudo-random problem}
\label{fig:randEx-tot}\vspace*{-0.1cm}
\end{figure}

Figures~\ref{fig:randEx} and \ref{fig:randEx3D} show an example of generated pseudo-random 2D and 3D instances with $n=20$, respectively.

\begin{algorithm}[h]
{
\caption{Pseudo-random trajectory generation} \label{alg1}\vspace{-3mm}
\begin{longtable}{lll}
Input      & $\hat p_i$, $\theta^{min}_i, \theta^{max}_i$, $\phi^{min}_i, \phi^{max}_i$, $v^{min}_i, v^{max}_i$ for all $i=1,\ldots,n$;  $n,n_c,p_c,max_c,maxTrials$\\[0.2cm]
			
Output     & $\hat V_i$: aircraft vectors of velocity, $i=1,\ldots,n$. \\[0.2cm]

\textbf{Step 0} & $isExplored[i]:=false$, for all $i=1,\ldots,n$, $nExplored:=0$, $totalConf:=0$.\\[0.1cm]
\textbf{Step 1}  & While ($nExplored < n$) do: \\[0.1cm]
& \textbf{1.0}  Take $i$ such that $isExplored[i]=false$.\\[0.1cm]
& \textbf{1.1}  Generate randomly $\hat v\in[v^{min}_i, v^{max}_i]$, $\hat \theta\in[\theta^{min}_i, \theta^{max}_i]$, and $\hat \phi\in[\phi^{min}_i, \phi^{max}_i]$.\\[0.1cm]
& \textbf{1.2}  $maxConf:=\min\{max_c,nExplored,n_c-totalConf\}$.\\[0.1cm]   
& \textbf{1.3}  $targetConf:=0$.\\[0.1cm]
& \textbf{1.4}  If $nExplored \neq 0$: \\[0.1cm]
& \hspace{0.8cm}If $\lceil \frac{n_c-totalConf}{(max_c+1)/2} \rceil\geq n-nExplored$ then $targetConf:=maxConf$.\\[0.1cm]
& \hspace{0.8cm}Else with probability $p_c$ do: \\[0.1cm]
& \hspace{1.3cm} $targetConf:=randInt(1,maxConf)$.\\[0.1cm]
& \textbf{1.5} $numTrials:=0$; $numTargetsTried:=0$; $targetNotReached:=true$ \\[0.1cm]
& \textbf{1.6} While ($targetNotReached$) do: \\[0.1cm]
& \hspace{0.8cm} While ($numTrials< maxTrials$ and $countConf(i,isExplored) \neq targetConf$) do: \\[0.1cm]
& \hspace{1.3cm} Generate randomly $\hat v\in[v^{min}_i, v^{max}_i]$, $\hat \theta\in[\theta^{min}_i, \theta^{max}_i]$, and $\hat \phi\in[\phi^{min}_i, \phi^{max}_i]$.\\[0.1cm]
& \hspace{1.3cm}  $numTrials++$.\\[0.1cm]
& \hspace{0.8cm} If ($countConf(i,isExplored)=targetConf$) do:\\ [0.1cm]
& \hspace{1.3cm} $targetNotReached:=false$\\[0.1cm]
& \hspace{0.8cm} Else do:\\ [0.1cm]
& \hspace{1.3cm} $numTargetsTried++$\\[0.1cm]
& \hspace{1.3cm} If $numTargetsTried\geq maxConf+1$ then $targetConf:=numTargetsTried$\\[0.1cm]
& \hspace{1.3cm} Else chose a non-tried $targetConf$ in $\{0,\ldots,maxConf\}$\\[0.1cm]
& \hspace{1.3cm} $numTrials:=0$\\[0.1cm]
& \textbf{1.7} $\hat{V}_i:=\hat  v(\cos\hat \theta \sin\hat \phi,\sin\hat \theta \sin\hat \phi,\cos\hat \phi)$ \\[0.1cm]
& \textbf{1.8} $totalConf:= totalConf+targetConf$ \\[0.1cm]
& \textbf{1.9} $nExplored++$ \\[0.1cm]
& \textbf{1.10} $isExplored[i]:=true$ \\[0.1cm]
& \textbf{1.11} If $nExplored \neq n$ \\[0.1cm]
& \hspace{0.7cm}$p_c=4\frac{n_c-totalConf}{(n-nExplored)*(1+max_c)}$
\end{longtable}
}
\end{algorithm}

\clearpage
\newpage

\section{Computational tests}\label{sec:compu}

\setcounter{table}{0}

\begin{table}[h]
    \begin{minipage}{0.48\linewidth}\hspace*{-0.1cm}
    \centering
    \begin{tabular}{r|rrrrr} 
      & \multicolumn{5}{c}{$den$, $n_c=den*\frac{n(n-1)}{2}$}  \\
        \hline
          $n$   & 0.05 & 0.10   & 0.15  & 0.20 & 0.25      \\
        \hline
        10	& 0	  & 3.75	& 0.71	& 3.89	& 6.82	   \\
        15	& 0	  & 1.00	& 0	    & 0.24	& 2.98	   \\
        20	& 0	  & 0.26	& 0.36	& 6.51	& 16.30	   \\
        25	& 0.33  &	0.08	& 2.61	& 9.58	& 20.27\\
        30	& 0.45  &	0.11	& 2.50	& 14.51	& 25.73\\
        35	& 0.42  &	0.21	& 2.42	& 11.47	& 23.27\\
        40	& 3.01  &	0.35	& 2.80	& 14.97	& 23.18\\
        45	& 2.25  &	0.33	& 1.71	& 11.67	& 24.65\\
        50	& 4.02  &	0.92	& 1.79	& 8.45	& 18.64\\
        75	& 10.05 &	2.45	& 8.13	& 20.50	& 26.35\\
        100	& 15.14 &	2.35	& 11.55	& 25.99	& 36.55\\
    \end{tabular}
 \subcaption{\textbf{2D} pseudo-random scenario}
  \label{mean-2D}
  \end{minipage}\hfill
    \begin{minipage}{0.4\linewidth}
    \centering
    \begin{tabular}{r|rrrrr}
      & \multicolumn{5}{c}{$den$, $n_c=den*\frac{n(n-1)}{2}$}  \\
        \hline
          $n$   & 0.05 & 0.10   & 0.15  & 0.20 & 0.25      \\
        \hline
        10 &	0	&	0	&	3.57	&	9.72	&	12.95	\\
        15 &	0	&	0.25	&	2.34	&	9.17	&	20.29	\\
        20 &	0.50	&	0.92	&	3.21	&	7.63	&	16.98	\\
        25 &	0.00	&	1.67	&	9.22	&	21.71	&	32.30	\\
        30 &	0.57	&	1.99	&	7.42	&	17.47	&	29.54	\\
        35 &	0.08	&	3.42	&	10.17	&	23.80	&	32.35	\\
        40 &	0	&	2.63	&	14.49	&	25.34	&	36.79	\\
        45 &	0	&	1.69	&	5.86	&	14.57	&	24.26	\\
        50 &	0.41	&	8.95	&	23.44	&	33.91	&	43.08	\\
        75 &	0.02	&	8.77	&	21.56	&	31.36	&	39.34	\\
        100 &	0.64	&	11.97	&	30.59	&	43.05	&	53.01	\\
    \end{tabular}
  \subcaption{\textbf{3D} pseudo-random scenario}
  \label{mean-3D}
  \end{minipage}
  \caption{Mean of the relative difference between number of obtained conflicts and $n_c$, over 40 pseudo-random instances.}
\end{table}

\begin{table}[ht]
\begin{minipage}{0.4\linewidth}\hspace*{-0.7cm}
\centering
     \begin{tabular}{r|rrrrr}
      & \multicolumn{5}{c}{$den$, $n_c=den*\frac{n(n-1)}{2}$}  \\
    \hline
        $n$ & 0.05 & 0.10   & 0.15  & 0.20 & 0.25      \\
    \hline
10&	0 (40)&	0 (34)&	0 (38)&	0 (26)&	0 (22)	\\
15&	0 (40)&	0 (36)&	0 (40)&	0 (39)&	0 (26)	\\
20&	0 (40)&	0 (39)&	0 (36)&	0 (15)&	0 (5)	    \\
25&	0 (38)&	0 (39)&	0 (29)&	0 (12)&	6.67 (1)	\\
30&	0 (36)&	0 (39)&	0 (25)&	0 (1)	&0.92 (1)	    \\
35&	0 (35)&	0 (36)&	0 (23)&	0 (7)	&1.34 (1)	    \\
40&	0 (22)&	0 (34)&	0 (22)&	0 (2)	&8.21 (1)	    \\
45&	0 (27)&	0 (33)&	0 (21)&	0 (6)	&0 (1)	        \\
50&	0 (23)&	0 (30)&	0 (28)&	0 (6)	&6.54 (1)	    \\
75&	0 (7)	&0 (22)&	0 (8)&	0 (7)&	0 (6)	                \\
100&	0.81 (1)&	0 (19)&	0 (7)&	0 (4)&	0.16 (1)	        \\
    \end{tabular}
    \subcaption{\textbf{2D} pseudo-random scenario}
    \label{min-2D}
\end{minipage}\hfill
\begin{minipage}{0.4\linewidth}\hspace*{-1cm}
\centering
    \begin{tabular}{r|rrrrr}
    & \multicolumn{5}{c}{$den$, $n_c=den*\frac{n(n-1)}{2}$}  \\
    \hline
    $n$ & 0.05 & 0.10   & 0.15  & 0.20 & 0.25      \\
    \hline
    10 &	0 (40)	&	0 (40)	&	0 (30)	&	0 (16)	  &	0 (18)	\\
    15 &	0 (40)	&	0 (39)	&	0 (27)	&	0 (17)	  &	0 (2)	\\
    20 &	0 (38)	&	0 (37)	&	0 (26)	&	0 (18)	  &	0 (8)	\\
    25 &	0 (40)	&	0 (29)	&	0 (12)	&	0 (5) 	  &	4 (1)	\\
    30 &	0 (39)	&	0 (31)	&	0 (22)	&	0 (10)	  &	0 (4)	\\
    35 &	0 (39)	&	0 (26)	&	0 (15)	&	0 (9)	  &	0 (5)	\\
    40 &	0 (40)	&	0 (27)	&	0 (10)	&	0 (5)     &	0 (1)	\\
    45 &	0 (40)	&	0 (31)	&	0 (25)	&	0 (12)	  &	0 (3)	\\
    50 &	0 (36)	&	0 (17)	&	0 (9)	&	0 (5)	  &	0 (1)	\\
    75 &	0 (39)	&	0 (20)	&	0 (14)	&	0 (9)	  &	1.44 (1)	\\
    100 &	0 (31)	&	0 (16)	&	0 (1)	&	18.08 (1) &	30.37(1)	\\
    \end{tabular}
    \subcaption{\textbf{3D} pseudo-random scenario}
  \label{min-3D}
    \end{minipage}
    \caption{Best (minimum) relative difference between number of obtained conflicts and $n_c$, over 40 pseudo-random instances. The number within parenthesis indicates at how many instances the minimum is attained.}
\end{table}

\begin{table}[ht]
    \begin{minipage}{0.48\linewidth}\hspace*{-0cm}
    \centering
    \begin{tabular}{r|rrrrr} 
      & \multicolumn{5}{c}{$den$, $n_c=den*\frac{n(n-1)}{2}$}  \\
        \hline
          $n$   & 0.05 & 0.10   & 0.15  & 0.20 & 0.25      \\
        \hline
        10 & 2.75 & 2.37 & 2.16 & 2.36 & 2.42 \\
        15 & 2.61 & 2.57 & 2.57 & 2.59 & 2.67 \\
        20 & 2.82 & 2.64 & 2.69 & 2.81 & 2.80 \\
        25 & 2.56 & 2.51 & 2.61 & 2.64 & 2.67 \\
        30 & 2.80 & 2.67 & 2.65 & 2.69 & 2.68 \\
        35 & 2.59 & 2.63 & 2.63 & 2.65 & 2.69 \\
        40 & 2.67 & 2.61 & 2.61 & 2.67 & 2.66 \\
        55 & 2.61 & 2.59 & 2.60 & 2.63 & 2.66 \\
        50 & 2.55 & 2.59 & 2.62 & 2.64 & 2.67 \\
        75 & 2.49 & 2.53 & 2.59 & 2.60 & 2.62 \\
        100 & 2.52 & 2.52 & 2.56 & 2.59 & 2.59 
    \end{tabular}
  \subcaption{\textbf{2D} pseudo-random scenario}
  \label{distance-2D}
  \end{minipage}\hfill
    \begin{minipage}{0.4\linewidth}
    \centering
    \begin{tabular}{r|rrrrr}
      & \multicolumn{5}{c}{$den$, $n_c=den*\frac{n(n-1)}{2}$}  \\
        \hline
          $n$   & 0.05 & 0.10   & 0.15  & 0.20 & 0.25      \\
        \hline
        10 &	3.62 &	3.39 &	3.47 &	3.37 &	3.36 \\
        15 &	3.55 &	3.59 &	3.45 &	3.37 &	3.44 \\
        20 & 	3.25 &	3.36 &	3.35 &	3.36 &	3.32 \\
        25 & 	3.42 &	3.34 &	3.39 &	3.36 &	3.39 \\
        30 & 	3.36 &	3.40 &	3.36 &	3.33 &	3.37 \\
        35 & 	3.37 &	3.32 &	3.35 &	3.31 &	3.29 \\
        40 & 	3.39 &	3.32 &	3.30 &	3.39 &	3.38 \\
        45 & 	3.42 &	3.36 &	3.37 &	3.33 &	3.33 \\
        50 & 	3.38 &	3.36 &	3.37 &	3.36 &	3.34 \\
        75 & 	3.35 &	3.32 &	3.28 &	3.28 &	3.26 \\
        100 & 	3.33 &	3.30 &	3.29 &	3.27 &	3.28
    \end{tabular}
  \subcaption{\textbf{3D} pseudo-random scenario}
  \label{distance-3D}
  \end{minipage}
  \caption{Mean of the average (over all the pairs in conflict) pairwise aircraft distance [NM] at the time of minimal separation, over 40 pseudo-random instances.}
\end{table}

\begin{table}[ht]
    \begin{minipage}{0.48\linewidth}\hspace*{-0cm}
    \centering
    \begin{tabular}{r|rrrrr} 
      & \multicolumn{5}{c}{$den$, $n_c=den*\frac{n(n-1)}{2}$}  \\
        \hline
         $n$   & 0.05 & 0.10   & 0.15  & 0.20 & 0.25      \\
        \hline
        10 &	1.35 &	1.21 &	1.54 &	14.80 &	2.83 \\
        15 &	2.68 &	2.39 &	3.75 &	4.12 &	2.05 \\
        20 & 	1.59 &	2.67 &	8.93 &	8.67 &	5.19 \\
        25 & 	2.85 &	5.28 &	4.23 &	4.78 &	4.09 \\
        30 & 	3.02 &	5.57 &	18.04 &	20.87 &	3.63 \\
        35 & 	3.02 &	2.79 &	4.98 &	8.40 &	2.89 \\
        40 & 	3.34 &	5.41 &	3.37 &	2.95 &	2.87 \\
        45 & 	4.01 &	5.85 &	8.15 &	4.38 &	4.46 \\
        50 & 	4.18 &	4.42 &	11.73 &	8.62 &	8.95 \\
        75 & 	4.58 &	5.42 &	6.72 &	5.29 &	8.11 \\
        100 & 	3.69 &	6.11 &	4.90 &	3.90 &	4.07
    \end{tabular}
  \subcaption{\textbf{2D} pseudo-random scenario}
  \label{duration-2D}
  \end{minipage}\hfill
    \begin{minipage}{0.4\linewidth}
    \centering
    \begin{tabular}{r|rrrrr}
      & \multicolumn{5}{c}{$den$, $n_c=den*\frac{n(n-1)}{2}$}  \\
        \hline
          $n$   & 0.05 & 0.10   & 0.15  & 0.20 & 0.25      \\
        \hline
         10 &	0.82 &	0.76 &	0.75 &	0.81 &	0.75 \\
        15 &	1.19 &	0.93 &	0.95 &	0.94 &	0.91 \\
        20 & 	1.12 &	0.84 &	0.96 &	0.80 &	0.89 \\
        25 & 	1.08 &	1.13 &	1.11 &	1.08 &	1.13 \\
        30 & 	0.92 &	1.03 &	0.99 &	0.93 &	0.90 \\
        35 & 	1.07 &	1.04 &	1.07 &	1.07 &	1.14 \\
        40 & 	1.20 &	1.22 &	1.22 &	1.21 &	1.22 \\
        45 & 	0.91 &	0.99 &	0.91 &	0.92 &	0.94 \\
        50 & 	1.20 &	1.21 &	1.11 &	1.17 &	1.32 \\
        75 & 	1.06 &	1.06 &	1.05 &	1.04 &	1.02 \\
        100 & 	1.17 &	1.14 &	1.13 &	1.14 &	1.15
    \end{tabular}
  \subcaption{\textbf{3D} pseudo-random scenario}
  \label{duration-3D}
  \end{minipage}
  \caption{Mean of the average (over all the pairs in conflict) pairwise duration [min] of conflict, over 40 pseudo-random instances.}
\end{table}

To demonstrate the efficacy of the pseudo-random generator in meeting the requested traffic congestion, we perform a computational study. For different numbers of aircraft $n$ and requested conflicts $n_c$, we vary some of the other options of the generator---namely, $max_c$ and the size of the air sector---and check whether the resulting instances feature the desired congestion. To further characterize the hardness of generated conflicts, we analyse two other indicators: the average pairwise distance at the time of minimal separation, and the average pairwise conflict duration (i.e., for how long the distance between aircraft is less then the given safety threshold of $5$ NM).

The experimental setting is the following. The number of aircraft is $n\in\{10,15,20,25,$ $30,35,40,45,50,75,100\}$. For a fixed $n$, we consider $n_c$ to correspond to different densities of conflicts. The density, which we denote by $den$ ($den\in[0,1]$) is the proportion with respect to the maximum number of possible conflicts, which is equal to $(n\cdot(n-1))/2$. The number of conflicts $n_c$ is the result of rounding $den\cdot(n\cdot(n-1))/2$ to its closest integer. 

We test several values of density, $den\in\{0.05,0.10,0.15,0.20,0.25\}$. The maximum number of conflicts per aircraft, $max_c$, is also tuned: $max_c=s+t$, with $s$ being the closest integer to $4\cdot(n_c/n)$ and $t\in\{1,2,3,4,5\}$. Finally, we change the size of the air sector, which is: 
\begin{itemize}
    \item for the 2D case, a $w\times w$ square, with $w\in\{125, 150, 175, 200, 225, 250, 275, 300\}$;
    \item for the 3D case, a $w\times w \times w$ cube, with $w\in\{50,60,70,80,90,100,125,150\}$.
\end{itemize}
The remaining options are set to their default values. This makes a total of 2200 running tests for the 2D case and 2200 running tests for the 3D case.

Tables~\ref{mean-2D} and \ref{mean-3D} show average results on the relative difference between $n_c$ and the number of conflicts featured by the generated instance, say $n'_c$. A complementary information is that, when $n'_c$ is different from $n_c$, in most cases $n'_c<n_c$. Tables \ref{min-2D} and \ref{min-3D} report the minimum of such relative differences that we obtain when $max_c$, and the size of the sector are varied. The number of configurations (among the 40 tried in total) for which this minimum value is attained is shown in parenthesis. If we use subscripts $i$ to denote a given instance with $n^i$ aircraft, $n_c^i$ requested conflicts and resulting ${n'}_c^i$ conflicts, Tables \ref{mean-2D} and \ref{mean-3D} show
$$\frac{1}{40}\sum_{i=1}^{40} 100\cdot \frac{|n_c^i-{n'_c}^i|}{n_c^i},$$
while Tables \ref{min-2D} and \ref{min-3D} show
$$\min_{i \in \{1,\dots,40\}}\left\{100\cdot\frac{|n_c^i-{n'_c}^i|}{n_c^i}\right\}.$$
For each entry of the tables, these mean and minimum values are calculated over 40 instances (which account for 8 possible sector sizes and 5 possible values of $max_c$).

Tables~\ref{distance-2D} and \ref{distance-3D} report the means of the minimal pairwise aircraft separation (in NM) averaged over all the pairs in conflict  for each tested instance. The minimum safety distance is set to $5$ NM. The tables show that in both 2D and 3D cases the minimal separation is around $3$ NM. The information on how long, in average, the aircraft pairs are in conflict can be found in Tables~\ref{duration-2D} and \ref{duration-3D}. 

We can notice that it is more difficult to generate the target number of conflicts when we take a third dimension into account, and also that the duration of the conflicts is greater for 2D instances. This comes as a natural consequence of problem dimensions. By adding an extra degree of freedom, the aircraft have an extra dimension to ``scape'' conflicts.
Despite this, Table \ref{mean-3D} shows values below 53\% of averaged relative difference between the number of obtained and requested conflicts. Moreover, from Table \ref{min-3D} we can notice that,  for most of the tested instances (51 out of 55), at least one combination of parameters exists to achieve the desired congestion. 

\section{Conclusion and future extensions}\label{sec:conclusion}
In this paper, we have presented an aircraft conflict resolution instances generator. The generator can produce different types of benchmarks: the well-known scenario-based instances, the randomly generated ones, and new pseudo-random configurations that we introduce in this work. The latter consists in
instances generated to meet a particular level of congestion, which is determined by the user through different input parameters. The algorithm used to generate these type of instances has been presented, and, to prove its efficacy in meeting the requested traffic congestion, several experiments have been conducted.
These bring light into the algorithm strengths and limitations. 

Our experiments demonstrate the utility of the proposed approach to generate instances meeting the desired level of congestion. The main factors influencing the accuracy of the method are the number of aircraft and the density of conflicts. The higher they are, the more difficult becomes to generate the target number of conflicts, i.e., more calls to the algorithm will be needed.
This reveals that there is room for improvement of the algorithm.
To this end, we identify two possible lines of amelioration, although we leave them as future work. One of them is to implement some kind of enhancement on top of the current algorithm. For instance, once a scenario is generated, we could go and try speed adjustments or small heading changes to generate more conflicts (if needed). A second line is to analyse the different random decisions of the algorithm and try to enhance them. This includes, for instance, the generation of the initial positions, and the order in which the vectors of velocity of the aircraft are generated. To this end, assessing the empirical distribution of conflicts among  aircraft could help to identify and correct possible imbalances of the current approach. In the same line, different distributions of the target number of conflicts other than the uniform distribution can be worth exploring.
A further new feature for future releases of the generator is guaranteeing some desired degree of complexity of the generated instances, not only with respect to the traffic density, but also to other indicators (e.g., conflict duration).

\paragraph{\textbf{Acknowledgments:}}{This work was partly funded by the European Union’s Horizon 2020 research and innovation programme under the Marie Sklodowska-Curie grant agreement n.\ 764759 ITN ``MINOA''.}

\bibliography{Cerulli_generator} 
\newpage
\appendix
\section{Parameters and notation}\label{appendix1}
\begin{table}[H]
\scalebox{0.8}{
    \centering
    \begin{tabular}{r|l}
        \rowcolor{light-gray}
        \textbf{Parameter} & \textbf{Definition} \\
        \hline
        $n$ & Number of aircraft \\
        \hline
        $\hat p_i = (\hat x_i, \hat y_i, \hat z_i)$ & Initial position of aircraft $i$ \\
        \hline
        $\hat V_i$ & Vector of velocity of aircraft $i$ \\
        \hline
        $\hat v_i$ & Speed of aircraft $i$ (i.e., modulus of $\hat V_i$)\\
        \hline
        $[v^{min}_i,v^{max}_i]$ & Feasible range of speed of aircraft $i$ \\
        \hline
        $\hat \theta_i$ & Heading angle with $x$-axis of $\hat V_i$ \\
        \hline
        $[\theta^{min}_i, \theta^{max}_i]$ & Feasible range of $\hat \theta_i$\\
		\hline
		$\hat \phi_i$ & Heading angle with $z$-axis of $\hat V_i$ \\
        \hline
		$[\phi^{min}_i, \phi^{max}_i]$ & Feasible range of $\hat \phi_i$ ($\phi^{min}_i=\phi^{max}_i=\pi/2 \quad \forall i=1,\ldots,n$ in the 2D case) \\
		\hline
		$D$ & Minimum safety distance \\
		\hline
		{$\alpha$} & Slope of slopping trails in the rhomboidal/grid scenario \\
		\hline
		HP & Number of horizontal planes in the polyhedral/cubic scenario \\
		\hline
		VP & Number of vertical planes in the polyhedral/cubic scenario \\
		\hline
		{$\alpha^h$} & Slope of slopping trails in the horizontal plane $h$ in the polyhedral/cubic scenario \\
		\hline
		$\beta^v$ & Slope of slopping trails in the vertical plane $v$ in the polyhedral/cubic scenario\\
		\hline
        {$m_x$} & Number of horizontal trails in the rhomboidal/grid scenario \\
        \hline
        {$m_x^h$} & Number of horizontal trails in the horizontal plane $h$ in the polyhedral/cubic scenario \\
        \hline
        {$m_y$} & Number of slopping trails in the rhomboidal/grid scenario \\
        \hline
        $m_y^h$ & Number of slopping trails in the horizontal plane $h$ in the polyhedral/cubic scenario \\
        \hline
        $m_z^v$ & Number of slopping trails in the vertical plane $v$ in the polyhedral/cubic scenario \\
        \hline
        {$n_x$} & Number of aircraft on each horizontal trail in the rhomboidal/grid scenario 
        \\
        \hline
        $n_x^h$ & Number of aircraft on each horizontal trail in the horizontal plane $h$ in the polyhedral/cubic scenario \\
        \hline
        {$n_y$} & Number of aircraft on each slopping trail in the rhomboidal/grid scenario\\
        \hline
        {$n_y^h$}& Number of aircraft on each slopping trail in the horizontal plane $h$ in the polyhedral/cubic scenario\\
        \hline
        $n_z^v$ & Number of aircraft on each slopping trail in the vertical plane $v$ in the polyhedral/cubic scenario \\
        \hline
        $n_c$ & Number of pairs of aircraft that are in conflict \\
        \hline
        $p_c$ & Probability that one aircraft is in conflict with at least another aircraft \\
        \hline 
        $max_c$ & Maximum number of aircraft that a fixed aircraft can be in conflict with
        
    \end{tabular}}
    \caption{List of parameters}
    \label{tab:parameters}
\end{table}

\end{document}